\documentclass{article}

\usepackage{PRIMEarxiv}

\usepackage[utf8]{inputenc} 
\usepackage[T1]{fontenc}    
\usepackage{hyperref}       
\usepackage{url}            
\usepackage{booktabs}       
\usepackage{amsfonts}       
\usepackage{nicefrac}       
\usepackage{microtype}      
\usepackage{lipsum}
\usepackage{graphicx}
\graphicspath{{media/}}     

\pdfoutput=1

\usepackage{bbold}
\usepackage{color}
\usepackage{graphics} 
\usepackage{needspace}

\input{mysymbol.sty}
\usepackage{amsthm}
\usepackage{cite}
\usepackage{amsmath}
\usepackage{mathtools}
\usepackage{multirow}
\usepackage{url}
\usepackage{graphics, times, amsfonts}
\usepackage{tikz, epic,eepic}
\usetikzlibrary{shapes,arrows}
\usepackage{color}
\usepackage{epstopdf}
\usepackage{epsfig, amsbsy}
\usepackage{lipsum}
\usepackage{latexsym}
\usepackage{amscd, verbatim}
\usepackage{multirow}
\usepackage{booktabs}

\usepackage{tikz}

\usepackage{algorithm}
\usepackage{algpseudocode}
\usepackage{subfiles}
\newtheorem{theorem}{Theorem}

\newtheorem{proposition}{Proposition}

\newtheorem{lemma}{Lemma}
\newtheorem{definition}{Definition}

{\itshape}{\rmfamily}
\newtheorem{assumption}{Assumption}
\newtheorem{remark}{Remark}

\title{\LARGE 
{The Lagrangian Method for Solving Constrained Markov Games}}
\author{ Soham Das, Santiago Paternain, Luiz F. O. Chamon, and Ceyhun Eksin
\thanks{Soham Das and Ceyhun Eksin are with the Department of  Industrial and Systems Engineering at Texas A\&M University, College Station, TX, USA (Email:{\tt\small \; soham.das@tamu.edu; eksinc@tamu.edu}).
Santiago Paternain is with the Department of Electrical, Computer and Systems Engineering at Rensselaer Polytechnic Institute, Troy, NY, USA (Email:{\tt\small \; paters@rpi.edu}).
Luiz F. O. Chamon is with the Center for Applied Mathematics at Ecole Polytechnique, Palaiseau, France (Email:{\tt\small \; luiz.chamon@polytechnique.edu}).
This work was supported by NSF ECCS-1953694, NSF CCF-2008855, and CAREER 2239410.
}}

\usepackage{graphicx}
\graphicspath{ {./images/} }
\usepackage{amssymb}

\begin{document}

\maketitle

{
\begin{abstract}
\textbf{We propose the concept of a Lagrangian game to solve constrained Markov games. Such games model scenarios where agents face cost constraints in addition to their individual rewards, that depend on both agent joint actions and the evolving environment state over time. Constrained Markov games form the formal mechanism behind safe multiagent reinforcement learning, providing a structured model for dynamic multiagent interactions in a multitude of settings, such as autonomous teams operating under local energy and time constraints, for example. We develop a primal-dual approach in which agents solve a Lagrangian game associated with the current Lagrange multiplier, simulate cost and reward trajectories over a fixed horizon, and update the multiplier using accrued experience. This update rule generates a new Lagrangian game, initiating the next iteration. Our key result consists in showing that the sequence of solutions to these Lagrangian games yields a nonstationary Nash solution for the original constrained Markov game.
}

\end{abstract}
}
\section{Introduction} \label{sec_intro}
\noindent Imagine a world where a collective of autonomous agents, each capable of independent decision-making, must interact in a shared environment to achieve their respective goals. These agents could be cooperative, such as a mixed fleet of autonomous and human operated taxis coordinating to reduce customer pick-up waiting times \cite{xie2023two}, or competitive, like buyers/traders in a stock market striving to maximize individual long term profits \cite{9410223}, \cite{cont2024dynamics}. In many real-world applications, the agents are not given explicit instructions on how to interact. Instead, they must learn through experience, adapting their actions based on feedback from the environment and the behavior of other agents. Over time, they may develop belief models about environment dynamics, optimal individual strategies, and even forms of implicit coordination \cite{aeonCommitmentCooperation} to meet their objectives. 

\begin{figure}
    \centering
    \includegraphics[width=0.35\linewidth]{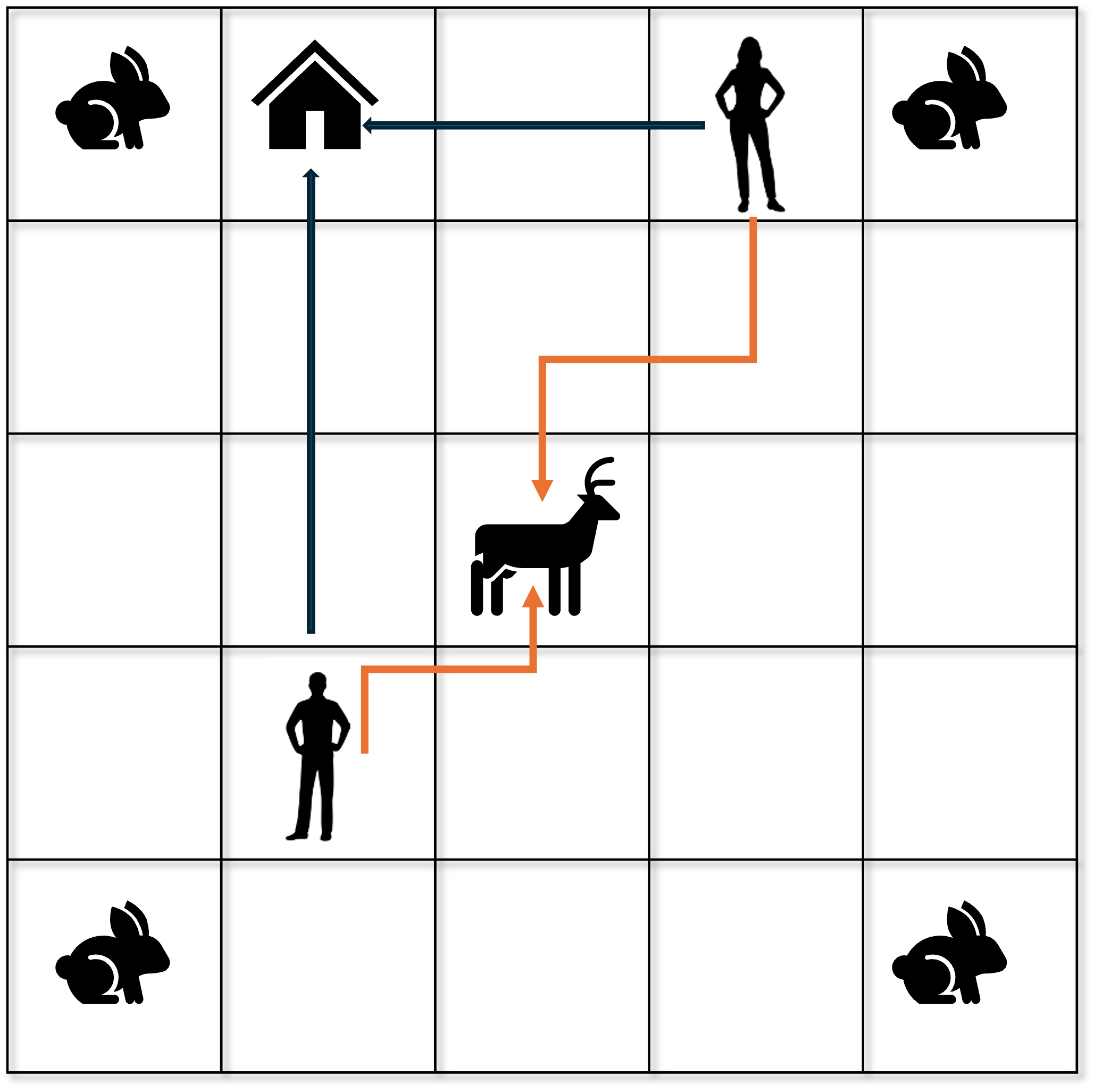}
    \caption{Agents can maximize payoffs by cooperating on hunting stag, but are required to spend a fraction of their hunting time at the resting station.}
    \label{fig:grid-world-hunting}
\end{figure}

For a concrete example, consider a grid-world hunting scenario, illustrated in Figure \ref{fig:grid-world-hunting}. We have two hunters on a foraging expedition, exploring the landscape in search of game. The grid contains locations where low-value targets, such as hares, can be hunted individually, while high-value targets, such as stags, require coordinated effort, that is, both hunters must arrive at the stag's location simultaneously to succeed.
Both individuals benefit most by cooperating to hunt a stag. However, each can also choose to hunt a hare alone. If one hunts a hare while the other pursues the stag, the hare-hunter gains a decent reward, while the stag-hunter gets nothing. If both hunt hares, they receive a lower payoff than if only one did, but at least they secure some reward.
Additionally, all hunters must allocate a fraction of their time to resting, ensuring long-term sustainability. The challenge, then, is twofold: agents must learn to coordinate effectively to maximize rewards while also complying with external resting constraints. A similar setting can be applied to coordinated surveillance, where autonomous agents must jointly monitor critical locations while managing energy consumption and resting requirements.

Multiagent Reinforcement Learning (MARL) has emerged as a powerful paradigm for modeling such complex, dynamic interactions. Through trial-and-error and learning from experience, agents can autonomously discover optimal strategies in multiagent environments.
MARL has demonstrated success in domains such as robotics, energy systems \cite{oroojlooy2023review}, and combinatorial games (such as Chess, Shogi \cite{silver2017mastering}), etc. However, as the grid-world example illustrates, real world multiagent interactions are rarely unconstrained. They are shaped by safety regulations, resource limitations, fairness considerations, etc. When these constraints are violated, systems may become unstable or unsafe, rendering standard MARL approaches inadequate.

Thus, in this paper, we shall advance the paradigm of \textit{safe}-Multiagent Reinforcement Learning (\textit{safe}-MARL), using the formal framework of constrained Markov games (CMGs). 
A Markov game (or stochastic game) \cite{shapley1953stochastic} involves repeated interactions among several participants when the environment state is dynamic and evolves in response to the actions of the agents in a stochastic fashion. Each player optimizes its own objective function while considering the actions of others. 
We consider repeated game play on an infinite time-horizon, with the system state evolving according to a transition probability kernel which satisfies the 
Markovian assumption. The agents take actions after each transition of the system with the goal to maximize their long run average payoffs, while ensuring compliance with potentially multiple other criteria. Such constrained problems arise in a wide range of real-world applications, such as in autonomous driving, where self-driving vehicles must balance speed and avoid collisions, or in sustainable energy management, where distributed energy grids must optimize power allocation while respecting environmental constraints.

\subsection{Related Works}

\noindent The study of Nash equilibria in constrained Markov games presents fundamental challenges due to the interplay between strategic decision-making, constraints, and the dynamic evolution of the game environment. 
Early foundational works have established the existence of Nash equilibria under various conditions \cite{nowak1985existence}, \cite{altman2000constrained}. Recently, Dufour et al. introduce the ARAT condition \cite{dufour2022stationary}, which facilitates equilibrium analysis. Their work on absorbing Markov games \cite{dufour2024nash} extends these results by addressing both constrained and unconstrained setups for the total expected reward value criterion, demonstrating the existence of stationary Nash in nonzero-sum games with measurable state spaces and compact metric action spaces. 

Despite theoretical advancements on existence of stationary Nash in constrained Markov games, most MARL research has focused on algorithm development for unconstrained settings, where agents optimize rewards without explicit feasibility constraints. This is, in part, because learning or computing approximate Nash (see Section \ref{sec_CMG} for related definitions) and other solution concepts in Markov games without constraints is itself a challenging problem. The approximate Nash problem is known to be PPAD-Complete in normal form games \cite{daskalakis2009complexity}, and therefore in Markov games.
In \cite{daskalakis2023complexity}, Daskalakis et al. prove that 
even computing stationary Markov Coarse Correlated Equilibria (CCE) is PPAD-Hard, utilizing the PCP for PPAD conjecture \cite{babichenko2015can}. CCE weakens the Nash Equilibrium solution concept by allowing for correlation between the strategies of agents. The authors circumvent the intractability result by instead providing a learning algorithm for \textit{nonstationary} Markov CCE. On the learning front, significant progress has been made on the polynomial computation and learning of approximate Nash in two player zero-sum Markov games. For example, independent policy gradient \cite{daskalakis2020independent}, V-Learning \cite{jin2021v}, self-play \cite{bai2020provable} and decentralized Q-learning \cite{sayin2021decentralized}  are effective for the zero-sum stochastic game model. 

When Markov games satisfy the structure that the difference in value corresponding to unilateral policy deviations by agents can be measured exactly by detecting changes in an associated potential function (Markov Potential Games or MPGs), policy gradient based approaches \cite{zhang2020global} become effective. Markov potential games naturally model (but are not limited to) identical-interest and cooperative scenarios. For examples, consider independent policy gradient \cite{ding2022independent} for large scale MPGs, multi-agent policy gradient \cite{leonardos2021global},  gradient-play \cite{zhang2024gradient} and two-time scale decentralized learning \cite{maheshwari2022independent} in MPGs. Beyond potential games, the convergence behaviors for policy gradient variants to special Nash policies (second order stationary) in general Markov games have been explored in \cite{giannou2022convergence}. Further, in \cite{mao2023provably}, Mao et al. develop a decentralized algorithm for provably learning approximate CCE in general sum Markov games. 


In the equilibirum learning literature for constrained Markov games, some recent results include the independent policy gradient based on proximal-point like policy updates \cite{jordan2024independent} and Nash coordinate ascent \cite{alatur2023provably} in constrained Markov Potential Games (CMPGs). 
In \cite{10531766}, Das et al. establish complexity guarantees for best-response dynamics for learning stationary Nash in constrained Markov games, when an upper bound for the maximum violation of the potential property (referred to as $\alpha$) is known for the game \cite{leonardos2021global}. Any Markov game can be equivalently interpreted as a Markov $\alpha$-potential game \cite{guo2023markov}, \cite{guo2024alpha}. In the present work, our proposed approach draws inspiration from constrained single-agent reinforcement learning (CRL) techniques and extends them to constrained Markov games. Specifically, we leverage a primal-dual framework comprising of a gradient descent step and policy optimization or learning to ensure feasibility in constrained dynamic multiagent environments.

In contrast to unconstrained MARL, CRL explicitly incorporates feasibility and safety requirements and is a well understood and successful paradigm. Constrained Markov decision processes (CMDPs) form the formal backbone for CRL problems. Several approaches have been proposed to enforce constraint satisfaction while optimizing rewards, such as linear programming \cite{6780601}, primal-dual \cite{vaswani2022near}, policy gradient \cite{liu2021policy} etc. Paternain et al. \cite{paternain2019constrained} established the foundational result that constrained reinforcement learning problems exhibit zero duality gap, explaining the success behind exact and approximate solutions for CRL in the dual domain \cite{chamon2022constrained}.
In \cite{calvo2023state}, Calvo-Fullana et al. introduced a novel state-augmented formulation for constrained reinforcement learning. The state space is augmented with the Lagrange multipliers and the primal-dual methods are reinterpreted as the portion of the dynamics that drives the multipliers evolution, 
generating systematic feasible solutions for CRL problems. Indeed, our game dynamics in this contribution are inspired by this result. Other approaches to CRL include adding probabilistic constraints for budget satisfaction, as developed in \cite{castellano2022reinforcement}, or probabilistic constraints for unsafe state visitation \cite{chen2024probabilistic}. Introducing constraints other than the standard expected discounted or ergodic time average reward introduces new opportunities for exploiting structure in the underlying CMDP---opportunities that remain unexplored as of yet in the Markov game literature.

\subsection{Our Contributions}
\noindent We introduce a novel framework for solving constrained Markov games by leveraging the formulation of a Lagrangian game. Given a Lagrange multiplier, the value functions of agents are augmented similarly to the construction of a Lagrangian in constrained optimization or CRL. While Lagrangian duality is a well-established concept in optimization, there has been no analogous paradigm for constrained Markov games. Our framework provides the first such formalization, enabling new algorithmic strategies for solving constrained multiagent decision-making problems.  A key result of our framework is that a trajectory of Lagrangian games can be solved sequentially, yielding a nonstationary constrained Nash equilibrium for the original constrained Markov game. 

The proposed algorithm follows a primal-dual approach. In the primal iteration, for a given value of the Lagrange multiplier, agents solve an unconstrained stationary Nash Equilibrium of the Lagrangian game using a suitable oracle. 
Thereby, the agents use the policy thus obtained to simulate reward trajectories for their value and cost functions over a finite time horizon. We call this simulation an epoch. Agents thereby utilize the experience obtained in the form of accrued constraint violations to update the Lagrange multipliers using stochastic gradient descent. The process repeats, refining the Lagrange multipliers and policies over multiple epochs.

Our analysis shows that the random cost trajectories obtained as a result of the continuous execution of the primal iteration in Algorithm 2 and the dual descent dynamics in Algorithm 1 results in episodes that are feasible with probability 1. Moreover, the sequence of  stationary Nash policies obtained as a result forms a nonstationary constrained Nash equilibrium for the constrained Markov game. Thus solving the harder constrained Markov game problem can be obtained by only solving unconstrained Lagrangian games, obtained from the original game by augmenting the value functions for the agents. We back the theoretical guarantees with practical solution results for an identical interest constrained Markov game. 

The proposed framework allows algorithm designers to solve constrained Markov games using existing solvers for unconstrained games. This significantly broadens the applicability of constrained game solutions without requiring novel constrained optimization techniques. Duality is not well-understood for constrained Markov games, as agents' payoff functions are often misaligned with each other. Unless the game is a fully cooperative game, there exists no joint objective function that agents seek to optimize, and therefore, no joint Lagrangian dual function for the agents. 
Additionally, it is now well-known that traditional strong duality does not hold even for CMPGs \cite{alatur2023provably}. Yet, our formulation of the Lagrangian game provides a systematic way to exploit duality arguments in general constrained Markov games, marking a first in the literature. The key idea is that in solving a sequence of Lagrangian games, where the Lagrange multipliers are updated following  Algorithm \ref{alg:dualdes}, we generate a policy sequence that is a nonstationary Nash solution to the original constrained Markov game. While our framework does not yield a stationary Nash solution, execution of the policy sequence generates cost-trajectories that are feasible almost surely. 

Complexity considerations for the proposed approach are outside the scope of the present article. Yet, it is clear that the complexity of our approach depends on the complexity of the underlying unconstrained Markov game solver. While our framework does not impose structural assumptions, specialized cases (such as Markov potential games \cite{leonardos2021global}, turn-based \cite{shah2020reinforcement},\cite{xie2020learning} or zero-sum Markov games \cite{sayin2021decentralized}) can lead to more efficient implementations, by harnessing structure in the associated Lagrangian games. Our framework, therefore, is flexible in that it does not prescribe a specific algorithm for solving the unconstrained Markov game. Moreover, it is compatible with both model-based and simulation-based approaches, making it applicable in both exact and reinforcement learning settings.


{\noindent{\it Notation:} The notations $\mbR$, $\mbN$ and $\mbZ$ represent, respectively, the set of real numbers, natural numbers and integers. $\mbR_{+}$ is the set of nonnegative reals and $\mbN_0$ is the set of natural numbers including 0. We use $\Delta(\ccalX)$ to denote the space of  probability distributions for any set $\ccalX$ (the probability simplex). We use brackets around an integer value $k$ to refer to the set $[k]:=\{1,2,\dots,k\}$. 
}

\section{Constrained Markov games} \label{sec_CMG}


\subsection{Game definition}
\noindent The constrained Markov game can be specified  by the tuple $\ccalG=(\ccalS,\ccalN,\{\ccalA_i, r_i\}_{i\in \ccalN}, P, \{c_{j},b_j\}_{j=1}^{m})$. 
Here $\ccalS$ is a finite set denoting the possible states of the system. We use $\ccalN=[N]$ to denote the set of $N\geq 2, N\in \mbN$ agents in the game. The set $\ccalA_i$ is a finite set representing the action space for agent $i\in \ccalN$ with elements $a_i\in \ccalA_i$. Denote by $\ccalA=\bigtimes_{i\in \ccalN}\ccalA_i$ the action space for all agents.
The function $r_i: \ccalS \times \ccalA \rightarrow \mbR$ is the individual reward function of agent $i\in \ccalN$. 
The game-play evolves over time $t\in \mbN_0$ (agents play the game repeatedly over time) and the global dynamic state $s\in\ccalS$ is driven by $P$, the transition probability kernel, i.e., $P(s'|s,a)$ is the probability of the state going from $s$ at time $t$ to state $s'$ at time $t+1$ when $a\in \ccalA$ is the action profile of the agents at time $t$. 


For each agent $i\in \ccalN$, we consider a stochastic stationary policy $\pi_i:\ccalS\rightarrow \Delta(\ccalA_i)$, where $\pi_i \in \Pi_i:=\Delta(\ccalA_i)^{\ccalS}$ determines a probability distribution over the actions of agent $i$ at each state $s\in \ccalS$. 
We define a policy function $\pi_{i}: \ccalS \times \mbN_0 \rightarrow \Delta(\ccalA_i)$ as stochastic nonstationary, when the probability distribution over the actions of agent $i$ are a function of the state $s \in \ccalS$ as well as the time index $t\in \mbN_0$ of the game evolution description. Thus $\pi_{i}\in \Pi^{\tau}_{i}:= \Delta(\ccalA_i)^{\ccalS\times \mbN_0}$, 
and $\pi_i(s,t_1)$ is not necessarily equal to $\pi_i(s,t_2)$ for $s\in \ccalS$ and $ t_1\neq t_2$. Equivalently, a stochastic nonstationary policy can be interpreted as a sequence of stochastic stationary policies, i.e. for $\pi_i\in \Pi^{\tau}_{i}$, we can write 
$\pi_i = \{\pi_{i}^{t}\}_{t\geq {0}}$ where $\pi_{i}^{t}:\ccalS\rightarrow \Delta (\ccalA_i)$ is stochastic stationary. 

Joint stochastic policies are defined similarly to individual player policies. A joint stationary stochastic policy is a mapping $\pi: \ccalS\rightarrow \Delta(\ccalA)$, whereas a joint stochastic nonstationary policy is a sequence $\pi=\{\pi^t\}_{t\geq 0}$, where each $\pi^t:\ccalS\rightarrow \Delta(\ccalA)$ is a stationary policy. Equivalently, $\pi:\ccalS\times \mbN_0 \rightarrow \Delta(\ccalA)$. Whenever the context is clear, we shall drop ``stochastic" and ``joint" from our terminology when referring to policies. We use $\Pi:= \Delta(\ccalA)^{\ccalS}$ and $\Pi^{\tau}:=\Delta(\ccalA)^{\mbN_0\times \ccalS}$ to denote the set of all stationary and nonstationary policies, respectively.
We say that a stationary policy $\pi$ is a \textit{product policy} when there are policies $\pi_i:\ccalS\rightarrow \Delta(\ccalA_i)$ such that $\pi(s)=\pi_1(s)\times...\times \pi_{N}(s)$ for all $s\in \ccalS$.
A nonstationary policy $\{\pi^t\}_{t\geq 0}$ is a product policy when each of its constituent policies $\pi^t:\ccalS\rightarrow \Delta(A)$ is a product policy. Given a stationary policy $\pi\in \Pi$, construct the following canonical mapping $\iota:\Pi\rightarrow \Pi^\tau$, where $\iota(\pi)=\{\pi\}_{t\geq 0}$. Then, $\iota(\Pi)\subseteq \Pi^\tau$, that is the set of stationary policies is a subset of the nonstationary policies under map $\iota$.

\noindent{\textit{The $-i$ notation}}: 
Given a stationary policy $\pi \in \Pi$ and a player $i\in \ccalN$, the mapping $\pi_{-i}:\ccalS \rightarrow \Delta(\ccalA_{-i})$ denotes the joint policy which at each state $s$ outputs the marginal distribution of $\pi(s)$ over $\ccalA_{-i}$, where $\ccalA_{-i}=\bigtimes_{j\in \ccalN, j\neq i} \ccalA_j$. Thus, $\pi_{-i}\in \Pi_{-i}:= \Delta(\ccalA_{-i})^{\ccalS}$.
For any $i\in \ccalN$ and stationary policies $\pi_i\in \Pi_i$ and $\pi_{-i}:\ccalS \rightarrow \Delta(\ccalA_{-i})$, the policy product $(\pi_i,\pi_{-i})$ refers to the policy which at each state $s$ samples an action profile according to product distribution $\pi_i(s)\times \pi_{-i}(s)$. When $\pi$ is a product policy, $\pi=(\pi_i,\pi_{-i})$.

For a joint nonstationary policy $\pi\in \Pi^\tau$, we denote $\pi_{-i}:=\{\pi_{-i}^{t}\}_{t\geq 0}$, where $\pi_{-i}^t:\ccalS\rightarrow \Delta(\ccalA_{-i})$.
For any $i\in \ccalN$ and nonstationary policies $\pi_i\in \Pi_{i}^{\tau}$ and $\pi_{-i}:\ccalS \times \mbN_0 \rightarrow \Delta(\ccalA_{-i})$, $(\pi_i,\pi_{-i})$ refers to the policy which at each state $s$ and time $t\geq 0$ samples an action profile according to product distribution $\pi_{i}^{t}(s)\times \pi_{-i}^{t}(s)$.

Given product policy $\pi=(\pi_i)_{i\in \ccalN}\in \Pi^{\tau}$, and some $i\in \ccalN$, we let $\pi_{-i}=(\pi_1,\pi_2,...,\pi_{i-1}, \pi_{i+1},...,\pi_{N}) \in \Pi^{\tau}_{-i} :=\bigtimes_{j\neq i, j\in \ccalN} \Pi_{j}^{\tau}$ denote the product policy of agents other than $i\in \ccalN$. 

\subsection{Value functions}
\noindent Consider a nonstationary policy $\pi \in \Pi^{\tau}$. Then the evolution of the stochastic game proceeds as follows. The initial state of the system is some $s^0\in\ccalS$ and at each step $t\in \mbN_0$, all players observe $s^t$ and draw actions $a^t$ according to policy $\pi^t$, where $\pi=\{\pi^t\}_{t\geq 0}$.
Then, the value function $V_i^{s}:\Pi^{\tau}\rightarrow \mbR$ gives the expected cumulative reward of agent $i\in \ccalN$ when $s^0=s$ and all agents draw their actions using policy $\pi$: 
\begin{align}
    V_i^{s}(\pi):=\liminf_{T\rightarrow \infty}\frac{1}{T}\mbE_{\pi}[\sum_{t=0}^{T} r_{i}(s^t, a^t) | s^0=s]. \label{eq_valuefn}
\end{align}
 In this work, we consider the expected time average reward (or infinite horizon average reward) criterion. Similarly, for a stationary policy $\pi\in \Pi$, the stationary value given by $V_i^{s}(\iota(\pi))$ returns the expected cumulative reward of agent $i\in \ccalN$, when $s^0=s$ and all agents draw their actions using policy $\pi$. {Whenever clear from the context, we drop the fixed initial state $s$ from the superscript of the value functions.}


\subsection{Constraints}

\noindent The constrained Markov game $\ccalG$ enforces that $m\in \mbN$ time-averaged system-level cost constraints be satisfied over repeated game play, for any nonstationary policy profile $\pi\in \Pi^\tau$ adopted by agents as follows
\begin{align}
    U_j^{s}(\pi):=\liminf_{T\rightarrow \infty}\frac{1}{T}\mbE_{\pi}\bigg[\sum_{t=0}^T c_{j}(s^t,a^t)|s^{0}=s\bigg] \geq b_j \textrm{ }\forall j\in [m],    \label{eq_constraints}
\end{align}
where $s\in \ccalS$ is the initial state, $s^t \in \ccalS$ and $a^t\in \ccalA$ denote the state and action profile at time $t\in \mbN_0$, $c_{j}:\ccalS\times \ccalA\rightarrow \mbR$ denotes the $j$-th cost function, and $b_j \in \mbR$ is the constant for constraint $j$. We define $b=(b_j)_{j\in [m]}$ and $U^s(\pi)=(U^s_{j}(\pi))_{j\in [m]}$. If agents adopt stationary policy $\pi\in \Pi$, the constrained Markov game enforces that $U^{s}_j(\iota(\pi))\geq b_j$ for all $j\in [m]$.

A policy $\pi$ is feasible, if it satisfies all  $m$-constraints. We allow the constraints to be \textit{coupled} across agents and time, that is they depend on the joint actions of all agents in the game, for all times. We use $\Pi^{c}\subseteq \Pi$ and $\Pi^{\tau,c} \subseteq \Pi^{\tau}$ to refer to the set of {stationary feasible policies} and {nonstationary feasible policies} respectively, in the game. Formally, $\Pi^{c}=\{\pi\in \Pi: U_j(\iota(\pi))\geq b_j \textrm{ } \forall j\in [m]\}$. We define $\Pi^{c}_i(\pi_{-i}):=\{\pi_i \in \Pi_i:(\pi_i,\pi_{-i})\in \Pi^{c}\}$ to refer to the set of feasible policies available to agent $i\in \ccalN$ when the remaining agents play $\pi_{-i}$. 
Similarly, we define $\Pi^{\tau,c}$ and $\Pi^{\tau,c}_{i}(\pi_{-i})$ as the nonstationary counterparts. 

\subsection{Solution concept: Approximate NE policy profile}

\begin{definition}(Nonstationary $\epsilon$-NE)\label{def_nonstatNE}
    A product nonstationary policy $\pi^{*}\in \Pi^{\tau,c}$ is a nonstationary $\epsilon$-NE of a constrained Markov game $\ccalG$ (or a constrained nonstationary $\epsilon$-NE), for some $\epsilon>0$, if for any $i\in \ccalN$ and $\pi_i\in \Pi_{i}^{\tau,c}(\pi_{-i}^{*})$ we have
    \begin{align}
        V^{s}_{i}(\pi^{*}) + \epsilon \geq V^{s}_{i}(\pi_{i},\pi^{*}_{-i})
    \end{align}
    {for all initial states} $s\in \ccalS$.
\end{definition}
 We refer to policy $\pi^{\star}_{i}$ an $\epsilon$-approximate \textit{best-response} to the joint policy of the other agents $\pi^{\star}_{-i}$. When $\epsilon=0$, we retrieve the standard NE definition.
In a similar vein, we define stationary $\epsilon$-NE as a solution concept for when agents adopt stationary policies.
\begin{definition}(Stationary $\epsilon$-NE)\label{def_statNE}
    A product stationary policy $\pi^{*}\in \Pi^{c}$ is a stationary $\epsilon$-NE of a constrained Markov game $\ccalG$ (or a constrained stationary $\epsilon$-NE), for some $\epsilon>0$, if for any $i\in \ccalN$ and $\pi_i\in \Pi_{i}^{c}(\pi_{-i}^{*})$ we have
    \begin{align}
        V^{s}_{i}(\iota(\pi^{*})) + \epsilon \geq V^{s}_{i}(\iota(\pi_{i},\pi^{*}_{-i}))
    \end{align}
    {for all initial states} $s\in \ccalS$.
\end{definition}

\noindent The NE is defined as a set of product policies of the players which satisfy simultaneously all the constraints and for which, in addition, no player can improve his payoff when unilaterally modifying his policy while still satisfying its own constraints. 
When the game $\ccalG$ does not have feasibility requirements, i.e. $\ccalG$ is an \textit{unconstrained} Markov game, then the Defs. \ref{def_nonstatNE}, \ref{def_statNE} correspond to the standard Defs. for nonstationary and stationary $\epsilon$-NE in Markov games. For notational convenience, whenever considering a stationary policy, we shall write $\iota(\pi)$ as $\pi$ keeping the $\iota(\cdot)$ map implicit. Thus, for instance, we shall write $V_i^{s}(\pi)$ instead of $V_i^{s}(\iota(\pi))$ when $\pi\in \Pi$ is stationary.


\section{The Lagrangian Method}
\subsection{Algorithmic framework}
 \begin{assumption}\label{asmp_oracle}
     Agents have access to an oracle that yields policies in the stationary NE of an unconstrained Markov game $\ccalG$. We use $NE(\ccalG)$ to represent the set of stationary NE policies for $\ccalG$.
 \end{assumption}
 We begin by stating that agents have a mechanism to access the stationary NE policies for a stochastic game that is unconstrained. Note that for a Markov game with finite number of actions and states, a stationary NE always exists. Methods such as value iteration, policy iteration and it's variants (such as natural policy gradient, actor-critic methods, etc.), and linear programming can be effective as methods for agents to compute/learn NE. Arguments related to time complexity of the various methods is not the current focus. Instead, given any approach to solve for NE for a unconstrained Markov game, we want to build on existing designs by satisfying system level constraints. 
Let us now define the important notion of the Lagrangian game. 
\begin{definition} (Lagrangian game) \label{def_lagrangainGame}
    Given the constrained Markov game $\ccalG$, the Lagrangian game $\ccalG(\lambda)$ is given by the tuple  $(\ccalS,\ccalN,\{\ccalA_i,r_i\}_{i\in \ccalN}, P, \{c_{j},b_j\}_{j=1}^{m}, \lambda)$, where $\lambda\in \mbR^{m}$ and the  value function $\ccalL_i: \Pi^{\tau} \times \mbR^{m}$ of agent $i\in \ccalN$ given by 
    \begin{align}
        \ccalL_i(\pi, \lambda):= V_i^{s}(\pi) + \lambda^{T}(U^{s}(\pi)-b). 
    \end{align}
\end{definition}
The algorithmic framework proceeds as follows.
Each episode of the repeated game play is divided into epochs of size $T_0$ time steps.
 Given constrained Markov game $\ccalG$, agents construct the Lagrangian game $\ccalG(\lambda_k)$ for epoch number $k \in \mbN_0$ (see algorithm \ref{alg:gamedyn}). 
 At epoch $k$, the agents access the oracle to solve $\ccalG(\lambda_k)$ and receive stationary NE policy $\pi^{k} \in NE(\ccalG(\lambda_k))$, following which they use $\pi^{k}$ to play the game $\ccalG(\lambda_k)$ (choose actions) for a horizon length of $T_0 \in \mbN$, the length of the epoch. Agents receive the roll-out rewards, which they use to construct $\lambda_{k+1}$ for the succeeding epoch as follows.
\begin{align} \label{eq_lambdaupdate}
    \lambda_{j,k+1}=\bigg[\lambda_{j,k}-\frac{\eta}{T_0}\sum_{t=kT_0}^{(k+1)T_{0}-1} (c_{j}(s^t,a^t)-b_j)\bigg]_{+}
\end{align}
for all $j\in [m]$. We refer to this as the dual descent update (see Algorithm \ref{alg:dualdes}). Here, $\lambda_k=(\lambda_{j,k})_{j\in [m]}$, and $\eta\in \mbR$ is the step-size of the update rule. At epoch $k+1$, agents construct Lagrangian game $\ccalG(\lambda_{k+1})$, and the process repeats itself. The successive updates to the Lagrange multiplier $\lambda_k$ lead to epoch dependent nonstationary policy $\pi^{\lambda} \in \Pi^{\tau}$, where $\pi^{\lambda}(s,t)=\pi^{k}(s)$ for $t\in [kT_0,(k+1)T_0)$, and all states $s\in \ccalS$. We shall abuse notation slightly to write $\pi^{\lambda}=\{\pi^{k}\}_{k\geq 0}$. The process is depicted in Algorithms \ref{alg:dualdes} and \ref{alg:gamedyn}.

\begin{algorithm} 
\caption{Dual Descent}
\label{alg:dualdes}
\begin{algorithmic}[1] 

\Require Constrained Markov game $\ccalG$, $\pi^{k} \in NE(\ccalG(\lambda_k))$, $\lambda_k \in \mbR^{m}_{+}$, $s^{kT_0} \in \ccalS$, $T_0 \in \mbN$, $\eta \in \mbR_{+}$.

\Function{dual-des} {$\pi^{k},s^{kT_0},\lambda_k$}
\State \textbf{Rollout} $T_0$ steps with actions $a^t \sim \pi^{k}(s^t)$ for $t\in [kT_0, (k+1)T_0-1]$
\ForAll{$j\in [m]$}
    \State $\lambda_{j,k+1}\gets\bigg[\lambda_{j,k}-\frac{\eta}{T_0}\sum\limits_{t=kT_0}^{(k+1)T_{0}-1} (c_{j}(s^t,a^t)-b_j)\bigg]_{+}$ 
\EndFor
\State \textbf{return} $(\lambda_{k+1}, s^{(k+1)T_0})$

\EndFunction
\end{algorithmic}
\end{algorithm}

\begin{algorithm} 
\caption{Game Dynamics}
\label{alg:gamedyn}
\begin{algorithmic}[1] 
\Require Constrained Markov game $\ccalG$, $s^{0}\in \ccalS$, $T_0 \in \mbN$, $\lambda_0 \in \mbR^{m}_{+}$.
\Function{play} {}
\ForAll {$\textrm{epoch }k=0,1,...K$}
    \State \textbf{Construct} $\ccalG(\lambda_k)$ 
    \State \textbf{Access} Oracle to solve $\ccalG(\lambda_k)$
    \State \textbf{Select} $\pi^{k} \in NE(\ccalG(\lambda_k))$
    \State $(\lambda_{k+1},s^{(k+1)T_0}) \gets$ \Call{DUAL-DES}{$\pi^{k},s^{kT_0},\lambda_k$}  
\EndFor
\EndFunction
\end{algorithmic}
\end{algorithm}

\subsection{Main result and implications}
\noindent To make claims for this algorithmic framework, we shall need the following assumptions.
\begin{assumption} \label{asmp_superfeasible}
    (Slater condition) For each $\pi\in\Pi$ and any player $i\in \ccalN$, there exists a $\pi^{\dagger}_i\in \Pi_i$ such that $U(\pi^{\dagger}_i,\pi_{-i})> b$. 
\end{assumption}
\begin{assumption} \label{asmp_boundedreward}
    (Bounded rewards) There exists $B$ such that $|c_{j}(s,a)-b_j|\leq B$ for all states $s\in \ccalS$, actions $a\in \ccalA$ and constraint $j\in [m]$. Moreover the reward functions $|r_{i}(s,a)|\leq R<\infty$ for all agents $i\in\ccalN$, for all states $s$ and actions $a$.
\end{assumption}
\begin{assumption} \label{asmp_unbiasedestimate}
    (Unbiased rollout) The accumulated reward that agents receive from the roll-out of stationary policy $\pi\in \Pi$ at any epoch $k$ is an unbiased estimate of $U_{j}(\pi)$, i.e.
    \begin{align}
        \mbE_{\pi} \left[ \frac{1}{T_0} \sum_{t=kT_0}^{(k+1)T_0-1} c_{j}(s^t,a^t)  \right] = U_{j}(\pi) 
    \end{align}
    for any constraint $j\in [m]$ and initial state. Similarly, we assume that 
    \begin{align}
        \mbE_{\pi} \left[ \frac{1}{T_0} \sum_{t=kT_0}^{(k+1)T_0-1} r_{i}(s^t,a^t)  \right] = V_i(\pi)
    \end{align}
    for all $i\in \ccalN$.
\end{assumption}
Assumption \ref{asmp_superfeasible} ensures that whatever policies other players use, a player $i\in \ccalN$ can find a policy to satisfy the constraints. This assumption is common in the constrained Markov game literature, see \cite{dufour2022stationary} for example.
Assumption \ref{asmp_boundedreward} ensures that the value functions are always finite.
In assumption \ref{asmp_unbiasedestimate}, we assume that the rollout reward at any epoch is an unbiased estimate of the value function under the executed policy at that epoch. This is a standard assumption in reinforcement learning, and we assume that $T_0$ is picked large enough for this assumption to be satisfied.
\begin{theorem} \label{thm_main}
    Under the above assumptions, the state action sequences $(s^t,a^t)$ generated by Algorithm \ref{alg:dualdes},\ref{alg:gamedyn} are feasible with probability $1$, i.e.
    \begin{align}
        \liminf_{T\rightarrow \infty} \frac{1}{T} \sum_{t=0}^{T-1} c_{j}(s^t,a^t) \geq b_j, \quad \forall j\in [m] \quad a.s.
    \end{align}
    Moreover, the sequence of policies forms a nonstationary $(\eta B^{2}/2)$-NE for the constrained Markov game $\ccalG$.
\end{theorem}
Thus, under assumptions \ref{asmp_oracle}-\ref{asmp_unbiasedestimate}, every realization of the policy sequence $\pi^{\lambda}=\{\pi^{k}\}_{k\geq 0}$ generated as a result of the continuous execution of the game dynamics in Algorithm \ref{alg:gamedyn} and the dual descent in Algorithm \ref{alg:dualdes} constitutes a constrained nonstationary approximate NE policy profile for the game $\ccalG$. 
Moreover, $\pi^{\lambda}$ has the structure that it is epoch dependent, as for epoch $k$, agents choose to play a stationary NE policy $\pi^{k}$ of the Lagrangian game $\ccalG(\lambda_k)$, for all time steps $t\in [kT_0,(k+1)T_0)$. The algorithmic framework is unique because we do not use \eqref{eq_lambdaupdate} as part of a learning process. Instead, we use \eqref{eq_lambdaupdate} as a online switching protocol that controls the change of the Lagrange multiplier, while the stationary NE policy $\pi^{k}$ corresponding to Lagrangian game $\ccalG(\lambda_k)$ is computed by the oracle. 

Note that the concept of \textit{duality} is not well established for general constrained Markov games.
Our results introduce a novel contribution by formalizing the notion of a Lagrangian game (Def. \ref{def_lagrangainGame}). Specifically, if agents can solve the unconstrained Lagrangian game, then Theorem \ref{thm_main} allows for the computation of a constrained nonstationary $\epsilon$-NE for the original constrained Markov game. 
The approach closely parallels primal dual optimization methods, where Algorithm \ref{alg:gamedyn} determines the ``primal" stationary policies and Algorithm \ref{alg:dualdes} updates the ``dual" Lagrange multipliers. 
\begin{remark}
    It is important to mention the caveat that the main result applies to the policy sequence and not to some limiting stationary policy in the sequence. In fact, we do not claim that the sequence of policies converges to some limiting stationary policy, which is feasible and approximate Nash.     
\end{remark}





\section{Analysis of the Lagrangian Method}
We start with definitions that will be used throughout the section. We define the probability space $(\Omega, \ccalF, \mbP)$ where the sample space is given by $\Omega = \ccalS \times \ccalA \times \mbR^{m}_{+}$, the event space $\ccalF$ corresponds to the Borel $\sigma$-algebra, and $\mbP:\ccalF\rightarrow [0,1]$ is the probability measure. We define a filtration $\{\ccalF_k\}_{k\geq 0}$ where $\ccalF_k \subset \ccalF$ is an increasing sequence of sigma algebras corresponding to the epoch $k$, i.e. if $k_2>k_1$ then $\ccalF_{k_1} \subset \ccalF_{k_2}$, with $\ccalF_0=\{\emptyset, \Omega\}$ and $\ccalF_\infty=\ccalF$. In particular, $\ccalF_k$ are such that $\lambda_k$ is $\ccalF_k$-measurable. We denote by $\{p(\lambda_k|\lambda_0)\}_{k\geq 0}$ the sequence of probability measures defined by the Lagrange multiplier dynamics in \eqref{eq_lambdaupdate}. We define the notion of tightness as follows.

 \begin{definition} \label{def_tightness}
    A sequence of probability measures $\{p(\lambda_k|\lambda_0)\}_{k\geq 0}$ are tight if for any $\delta>0$, there exists a set $\ccalK_\delta$ such that for each $k\geq 0$, we have $\mbP[\lambda_k\in \ccalK_{\delta}]>1-\delta$.
\end{definition}
Thus for any given small number $\delta$, there exists a compact set that captures nearly all the probability mass of each measure in the sequence.

\subsection{Proof of tightness}
In this section, we establish the tightness of the sequence of measures defined in \eqref{eq_lambdaupdate}, which forms the foundation for the results in the next subsections, which in turn enable Theorem \ref{thm_main}. 

\begin{definition}\label{def_gendual}
     We define the generalized dual $d_i: \mbR^{m}_{+} \times \Pi_{-i} \rightarrow \mbR$  for agent $i$ in constrained Markov game $\ccalG$ as follows
    \begin{align} \label{eq_dik1}
        d_i(\lambda,\pi_{-i})= \sup_{z_i \in \Pi_i} [V_i(z_i,\pi_{-i}) + \lambda^{T}(U(z_i,\pi_{-i})-b)]
    \end{align}
\end{definition}
The generalized dual represents the best that agent $i$ can achieve with respect to its Lagrangian value function, when the considered Lagrange multiplier is $\lambda\in\mbR^{m}_{+}$ and the policies of the other agents are fixed to the stationary policy $\pi_{-i}\in\Pi_{-i}$.
The two arguments for the generalized dual are the Lagrange multiplier and the stationary policy $\pi_{-i}$.

The following Lemma regarding generalized duals will be useful in establishing tightness. Variations of this result are often referred to as Danskin's Theorem in the literature. 

\begin{lemma} (Danskin's Theorem) \label{lem_danskin}
    For epoch $k$ in the execution of Algorithm \ref{alg:gamedyn}, with $\pi^k \in NE(\ccalG(\lambda_k))$ we have
    \begin{align}
        d_i(\lambda^{+},\pi^{k}_{-i})-d_i(\lambda_k,\pi^{k}_{-i}) \geq (\lambda^{+}-\lambda_k)^{T}(U(\pi^k)-b)\nonumber
    \end{align}
    where $\lambda^{+} \in \mbR^{m}_{+}$ is some reference Lagrange multiplier, and $i\in \ccalN$.
\end{lemma}
\begin{proof}
    At epoch $k$, agents other than $i\in \ccalN$ play policy $\pi^k_{-i}$. Define $L:\Pi_i\times \mbR^{m}_{+}\rightarrow \mbR$ as $L(z_i,\lambda)=V_i(z_i,\pi^{k}_{-i})+\lambda^{T}(U(z_i,\pi^{k}_{-i})-b)$. Then, using Def. \ref{def_gendual}, we have $d_i(\lambda,\pi^{k}_{-i})=\sup_{z_i \in \Pi_i} L(z_i,\lambda)$. Since $\pi^k\in NE(\ccalG(\lambda_k))$, $\pi^{k}_{i}$ lies in the best-response strategy set of agent $i$ against $\pi^{k}_{-i}$ for optimizing objective $L(\cdot,\lambda_k)$. Thus, 
    \begin{align}\label{eq_L1}        L(\pi^{k}_{i},\lambda_k)&=d_i(\lambda_k,\pi^{k}_{-i}) \nonumber \\
        \textrm{or, } d_i(\lambda^{+},\pi^{k}_{-i})&=d_i(\lambda^{+},\pi^{k}_{-i})+d_i(\lambda_k,\pi^{k}_{-i})-L(\pi^{k}_{i},\lambda_k)
    \end{align}
    where we add $d_i(\lambda^{+},\pi^{k}_{-i})$ to both sides on the second line in \eqref{eq_L1}. Also note that 
    \begin{align} \label{eq_L2}
        d_i(\lambda^{+},\pi^{k}_{-i}) \geq L(\pi^{k}_{i}, \lambda^{+})
    \end{align}
    as $\pi^{k}_i$ is a potentially suboptimal response for agent $i$ against objective $L(\cdot,\lambda^{+})$. As a result, plugging in \eqref{eq_L2} into \eqref{eq_L1}, we have,
    \begin{align}\label{eq_L3}
        d_i(\lambda^{+},\pi^{k}_{-i}) \geq L(\pi^k_{i},\lambda^{+}) + d_i(\lambda_k,\pi^{k}_{-i}) - L(\pi^k_{i},\lambda_k)
    \end{align}
    Then, rearranging terms and substituting for $L(\cdot,\cdot)$ we get the desired result.
\end{proof}


Let agents engage in game-play using nonstationary policy $\pi^{\lambda}$, where for epoch $k$, agents play $\pi^{k}$, a NE-policy for Lagrangian game $\ccalG(\lambda_k)$. Consider the following set 
\begin{align} \label{eq_Dk1}
\ccalD^{k}_{i}:=\{\lambda\in \mbR^{m}_{+}|d_i(\lambda,\pi^{k}_{-i})-d_i(\lambda^{+},\pi^{k}_{-i})\leq \eta B^2/2 \}    
\end{align}
Here $\lambda^{+}\in \mbR^{m}_{+}$ is some reference lambda which has finite norm, and $k$ is an arbitrary fixed epoch number.   

The set $\ccalD^{k}_i$ comprises of all $\lambda$ `close enough' to $\lambda^{+}$ in terms of the evaluation of the generalized dual for agent $i \in \ccalN$, for chosen right hand side $\eta B^2/2$.
Then, the following holds.
\begin{lemma}
    The set $\ccalD^{k}_i$ is contained in a compact set.
\end{lemma}
\begin{proof}
    By Slater's conditions (assumption \ref{asmp_superfeasible}), there exists feasible  policy $\pi^{\dagger}_i$ with respect to $\pi^{k}_{-i}$, which leads to
    \begin{align}
        d_i(\lambda,\pi^{k}_{-i}) \geq V_i (\pi^{\dagger}_i,\pi^{k}_{-i}) + \lambda^{T} (U(\pi^{\dagger}_i,\pi^{k}_{-i})-b) 
    \end{align}
    This implies, using $U(\pi^{\dagger}_i,\pi^{k}_{-i})-b\geq C>0$,
    \begin{align}
       d_i(\lambda,\pi^{k}_{-i})-d_i(\lambda^{+},\pi^{k}_{-i}) \\ \nonumber \geq  V_i (\pi^{\dagger}_i,\pi^{k}_{-i})+ \lambda^{T}C - d_i(\lambda^{+},\pi^{k}_{-i})
    \end{align}
    Now, for all $\lambda\in \ccalD^{k}_{i}$, we have $d_i(\lambda,\pi^{k}_{-i})-d_i(\lambda^{+},\pi^{k}_{-i})\leq \eta B^2/2$. Thus, we get that 
    \begin{align} \label{eq_Dkbound1}
    V_i(\pi^{\dagger}_i,\pi^{k}_{-i})+ \lambda^{T}C - d_i(\lambda^{+},\pi^{k}_{-i}) \leq \eta B^2/2.    
    \end{align}
    Let $C_{\star}=\min_{j\in [m]}C_j$, where $C=(C_j)_{j\in [m]}$. Since $\lambda\in \mbR^{m}_{+}$, and $C>0$, we have $\|\lambda\|_{1} C_{\star} \leq \lambda^{T}C$, where $\|\lambda\|_{1}=\sum_{j=1}^{m} \lambda_j$. Thus, using \eqref{eq_Dkbound1}, we get,
    \begin{align}
        \|\lambda\|_{1} C_{\star} \leq \eta B^2/2 + d_i(\lambda^{+},\pi^{k}_{-i}) -V_i(\pi^{\dagger}_i,\pi^{k}_{-i})
    \end{align}
    Since the reward functions $r_i$ are finite, we get that $\|\lambda\|_1$ is bounded, for all $\lambda\in \ccalD^{k}_{i}$. Then, $\ccalD_{i}^{k}$ is contained in the ball $\ccalD$ defined as follows.
    \begin{align}
        \ccalD=\{\lambda\in \mbR^{m}_{+}: \|\lambda\|_1 \leq \frac{\eta B^2/2+d_i(\lambda^{+},\pi^{k}_{-i})-V_i(\pi^{\dagger}_{i},\pi^{k}_{-i})}{C_{\star}}\}
    \end{align}
    Since the closed ball is compact, we have our desired result.
\end{proof}
\begin{lemma} \label{lem_tightness}
    Under the assumptions of Theorem \ref{thm_main}, the sequence of probability measures $\{p(\lambda_k|\lambda_0)\}_{k\geq 0}$ is tight according to definition \ref{def_tightness}.
\end{lemma}
\begin{proof}
    Let $\ccalD^{k}_{i}$ be the set previously defined, and $\bar\ccalD^{k}_{i}$ it's complement. By virtue of the previous Lemma, the set $\ccalD^{k}_{i}$ is contained in a compact set. Define $\ccalD^k:=\cup_{i\in \ccalN}\ccalD_i^{k}$. Then $\ccalD^k$ is also contained in a compact set.
    Let us denote by $\ccalB_{\eta B}$ a closed ball centered at $0$ and radius $\eta B$ and by $\oplus$ the Minkowski addition. 
    Given knowledge of sigma algebra $\ccalF_k$, consider the sequence of sets $(\ccalD^q)_{q\leq k}$. If $\lambda^k \in \ccalD^{k}$, then $\lambda_{k+1} \in \ccalD^{k} \oplus \ccalB_{\eta B}$, because the estimate of the gradient in \eqref{eq_lambdaupdate} is bounded by $B$. Since $\ccalB_{\eta B}$ is compact, we have that 
     $\ccalD^k \oplus \ccalB_{\eta B}$ is contained in a compact set. This means that the first iterate of the sequence $\{\lambda_k\}_{k\geq 0}$, say $\lambda_l$, that is in the set $\ccalD^l$, has bounded norm with probability $1$.

    Without loss of generality, assume that the sequence $\{\lambda_k\}_{k\geq 0}$ is such that $\lambda_0 \in \bar\ccalD^0$. Let us define the following stopping time $K_0=\min_{k\geq 0}\{\lambda_k \in \ccalD^k\}$. We define the sequence $\{\lambda_{k\wedge K_0}\}$, where $\wedge$ denotes the minimum between $k$ and $K_0$. Note that the previous discussion allows us to exclusively consider indices $k\leq K_0$. That is, $\lambda_k \in \ccalC$ for all $k\geq K_0$, where $\ccalC$ is some compact set.

    We denote by $\|\lambda-\lambda^{+}\|^{2}$ the distance from $\lambda$ to the reference $\lambda^{+}$, where $\|\cdot\|$ is the norm induced by the inner product. Using the nonexpansiveness of the norm, and substituting $\lambda_{(k+1)\wedge K_0}$ by its dynamics in \eqref{eq_lambdaupdate}, it follows that 
    \begin{align} \label{eq_tight1}
        \|\lambda_{(k+1)\wedge K_0} - \lambda^{+}\|^2 \leq \|\lambda_{k\wedge K_0}-\eta g_k -\lambda^{+}\|^{2} \nonumber \\
        = \|\lambda_{k\wedge K_0}-\lambda^{+}\|^{2} - 2\eta (\lambda_{k\wedge K_0}-\lambda^{+})^{T}g_k + \eta^2\|g_k\|^{2}
    \end{align}
    where the vector $g_k=\sum_{t=kT_0}^{(k+1)T_0-1} (c(s^t,a^t)-b)/T_0 $ is the time average of the difference between the constraint payoff rollout and the constraint right hand side. Next, consider the conditional expectation of \eqref{eq_tight1} with respect to the sigma algebra $\ccalF_k$.
    \begin{align} \label{eq_tight2}
        \mbE\left[ \|\lambda_{(k+1)\wedge K_0} - \lambda^{+}\|^2 | \ccalF_k \right] \leq \|\lambda_{k\wedge K_0} - \lambda^{+}\|^2 \nonumber \\ -2\eta (\lambda_{k\wedge K_0}-\lambda^{+})^{T}\mbE \left[g_k|\ccalF_k\right] +\eta^2 \mbE\left[\|g_k\|^2\right | \ccalF_k] .
    \end{align}
    Here we utilize the fact that $\lambda_{(k\wedge K_0)}$ is $\ccalF_k$-measurable. Using assumption \ref{asmp_boundedreward}, it follows that $\mbE\left[\|g_k\|^2\right | \ccalF_k] \leq B^2$. Plugging into \eqref{eq_tight2}, we get
    \begin{align}\label{eq_tight3}
        \mbE\left[ \|\lambda_{(k+1)\wedge K_0} - \lambda^{+}\|^2 | \ccalF_k \right] \leq \|\lambda_{k\wedge K_0} - \lambda^{+}\|^2 \nonumber \\ -2\eta (\lambda_{k\wedge K_0}-\lambda^{+})^{T}\mbE \left[g_k|\ccalF_k\right] +\eta^2 B^2.
    \end{align}
    Using assumption \ref{asmp_unbiasedestimate}, we have that $\mbE \left[g_k|\ccalF_k\right]=U(\pi^{k})-b$ for all $\pi^{k} \in NE(\ccalG(\lambda_k))$. 
    Therefore the second term on the RHS of \eqref{eq_tight3} is equivalent to $-2\eta(\lambda_{k\wedge K_0}-\lambda^{+})^{T}(U(\pi^{k})-b)$. By Danskin's Theorem (Lemma \ref{lem_danskin}), it follows that 
    $(\lambda^{+}-\lambda_{k\wedge K_0})^{T}\mbE\left[g_k|\ccalF_k\right] \leq d_i(\lambda^{+},\pi^{k}_{-i})-d_i(\lambda_{k\wedge K_0},\pi^{k}_{-i})$. Substituting in the previous inequality, we get
    \begin{align} \label{eq_tight4}
       \mbE\left[ \|\lambda_{(k+1)\wedge K_0} - \lambda^{+}\|^2 | \ccalF_k \right] \leq \|\lambda_{k\wedge K_0} - \lambda^{+}\|^2 \nonumber \\ -2\eta (d_i(\lambda_{k\wedge K_0},\pi^{k}_{-i}) - d_i(\lambda^{+},\pi^{k}_{-i})) +\eta^2 B^2. 
    \end{align}
    Since $\lambda_{k\wedge K_0}$ is outside $\ccalD^k$, it is outside $\ccalD_i^k$ for all agents $i\in \ccalN$. So we have that $(d_i(\lambda_{k\wedge K_0},\pi^{k}_{-i}) - d_i(\lambda^{+},\pi^{k}_{-i})) > \eta B^2 / 2$, which further implies that $-2\eta (d_i(\lambda_{k\wedge K_0},\pi^{k}_{-i}) - d_i(\lambda^{+},\pi^{k}_{-i})) +\eta^2 B^2 < 0$ in \eqref{eq_tight4}. Thus, we get that 
    \begin{align} \label{eq_tight5}
        \mbE\left[ \|\lambda_{(k+1)\wedge K_0} - \lambda^{+}\|^2 | \ccalF_k \right] \leq \|\lambda_{k\wedge K_0} - \lambda^{+}\|^2.
    \end{align}
    Applying this equation recursively and taking expectation, we get that 
    \begin{align} \label{eq_tight5B}
        \mbE\left[ \|\lambda_{(k+1)\wedge K_0} - \lambda^{+}\|^2 \right] \leq \|\lambda_{0} - \lambda^{+}\|^2.
    \end{align}

    Having established \eqref{eq_tight5B}, we construct set $\ccalK_\epsilon$, for any $\epsilon>0$ as follows
    \begin{align}
        \ccalK_\epsilon =\{\lambda \in \mbR^{m}_{+}: \|\lambda-\lambda^{+}\|^{2} \leq \|\lambda_0-\lambda^{+}\|^{2}/\epsilon\}.
    \end{align}
    The ball $\ccalK_\epsilon$ is closed by definition. Moreover, since $\lambda^{+}$ and $\lambda_0$ have finite norm, and the set $\ccalK_\epsilon$ consists of such $\lambda$ which are at most a fixed distance away from $\lambda^{+}$, $\ccalK_\epsilon$ is bounded. Thus $\ccalK_\epsilon$ is compact. Hence, using Markov inequality, we get that 
    \begin{align} \label{eq_tight6}
        \mbP\left[\lambda_{k\wedge K_0} \in \ccalK_\epsilon\right] = \mbP\left[ \|\lambda_{k\wedge K_0}-\lambda^{+}\|^{2} \leq \|\lambda_0-\lambda^{+}\|^{2}/\epsilon\right] \nonumber \\
        > 1-\epsilon \frac{\mbE\left[\|\lambda_{k\wedge K_0}-\lambda^{+}\|^{2}\right]}{\|\lambda_0-\lambda^{+}\|^{2}} \geq 1-\epsilon
    \end{align}
    where the last inequality follows from \eqref{eq_tight5B}. Thus, our proof is complete.
\end{proof}
\subsection{Feasibility guarantees assuming tightness}
\begin{proposition} \label{prop_feasibility1}
    If the sequence $\{p(\lambda_k|\lambda_0)\}_{k\geq 0}$ is tight, the state-action trajectories are feasible with probability 1, i.e.,
    \begin{align}
        \liminf_{T\rightarrow \infty}\frac{1}{T}\bigg[\sum_{t=0}^T c_{j}(s^t,a^t)|s^{0}=s\bigg] \geq b_j
    \end{align}
    holds almost surely for all $j\in [m]$ and for any arbitrary initial state $s$.
\end{proposition}
\begin{proof}
    For the sake of contradiction, assume there exists $\beta \in (0,1]$ and $\epsilon>0$ such that for some $j\in [m]$ we have 
    \begin{align} \label{eq_feasible1}
        \mbP\left[ \liminf\limits_{T\rightarrow \infty} \frac{1}{T} \sum_{t=0}^{T-1} c_{j}(s^t,a^t) \leq b_j -\epsilon\right]= \beta.
    \end{align}
    We first lower bound the norm of the $k+1$ iterate of the dual multiplier $\lambda_{j,k+1}$ whose evolution happens according to \eqref{eq_lambdaupdate}, as follows. See that by virtue of the projection onto the positive orthant in \eqref{eq_lambdaupdate}, we get 
    \begin{align} \label{eq_feasible2}
        \lambda_{j,k+1} \geq \lambda_{j,k} - \frac{\eta}{T_0}\sum_{t=kT_0}^{(k+1)T_0-1}\left(c_{j}(s^t,a^t) -b_j\right).
    \end{align}
    Applying \eqref{eq_feasible2} recursively, we get 
    \begin{align}\label{eq_feasible2-1}
        \lambda_{j,k+1} \geq \lambda_{j,0} - \frac{\eta}{T_0}\sum_{t=0}^{(k+1)T_0-1}\left(c_{j}(s^t,a^t) -b_j\right).
    \end{align}
    By the definition of the $l_1$ norm, one has that $\|\lambda_0\|_1 \geq \lambda_j$ for all $j\in [m]$ and $\lambda\in \mbR^{m}_{+}$. Hence, it holds that 
    \begin{align} \label{eq_feasible3}
        \limsup\limits_{k\rightarrow \infty} \|\lambda_{k+1}\|_1 \geq \limsup\limits_{k\rightarrow \infty} \lambda_{j,k+1}.
    \end{align}
    Using \eqref{eq_feasible2-1}, we obtain
    \begin{align} \label{eq_feasible3-1}
        &\limsup_{k\rightarrow \infty} \|\lambda_{k+1}\|_1 \\ \nonumber &\geq \lambda_{j,0} + \limsup\limits_{k\rightarrow \infty} \left[- \frac{\eta}{T_0}\sum_{t=0}^{(k+1)T_0-1}\left(c_{j}(s^t,a^t) -b_j\right)\right].
    \end{align}
    Now, 
    \begin{align}
        &\limsup\limits_{k\rightarrow \infty} \left[- \frac{\eta}{T_0}\sum_{t=0}^{(k+1)T_0-1}\left(c_{j}(s^t,a^t) -b_j\right)\right] \nonumber \\ &= -\eta\liminf\limits_{k\rightarrow \infty}(k+1)\frac{1}{(k+1)T_0}\left[\sum_{t=0}^{(k+1)T_0-1}\left(c_{j}(s^t,a^t) -b_j\right)\right].
    \end{align}
    Using \eqref{eq_feasible1}, we have that with probability $\beta$, 
    \begin{align}
        \liminf\limits_{k\rightarrow \infty}\frac{1}{(k+1)T_0}\left[\sum_{t=0}^{(k+1)T_0-1}\left(c_{j}(s^t,a^t) -b_j\right)\right] \leq -\epsilon,
    \end{align} 
    which would imply 
    \begin{align}
       \limsup\limits_{k\rightarrow \infty} \left[- \frac{\eta}{T_0}\sum_{t=0}^{(k+1)T_0-1}\left(c_{j}(s^t,a^t) -b_j\right)\right] \nonumber \\ \geq \eta \liminf\limits_{k\rightarrow \infty} (k+1)\epsilon = \infty. 
    \end{align}
    Thus we have shown that with probability $\beta$, the norm of the $(k+1)$ iterate of the Lagrange multiplier $\lambda_{k+1}$ is lower-bounded at infinity. Hence for any compact set $\ccalK$, there exists a constant $l_\ccalK\geq 0$ and a subsequence $\{k_l\}$ such that $\mbP(\lambda_{k_l}\in \ccalK)=1-\beta$ for all $l>l_\ccalK$. This further implies that for $\delta\in(0,\beta]$, there does not exist a compact set $\ccalK_\delta$ for which $\mbP(\lambda_k \in \ccalK_\delta) > 1-\delta$ for all iterates $k$, which contradicts our assumption that the sequence of probability measures $\{p(\lambda_k|\lambda_0)\}_{k\geq 0}$ is tight. This completes the proof. 
\end{proof}
\subsection{Optimality guarantees assuming tightness}
\begin{lemma} \label{lem_optimality1}
    Under the assumptions of Theorem \ref{thm_main}, when agents execute policy sequence  $\{\pi^{k}\}_{k\geq 0}$, ergodic complementary slackness holds, i.e.
    \begin{align}
        \limsup\limits_{K\rightarrow \infty} \frac{1}{KT_0} \sum_{k=0}^{K-1} \mbE \left[\lambda_k^{T} \sum_{t=kT_0}^{(k+1)T_0-1}(c(s^t,a^t)-b)\right] \leq \eta B^2 / 2
    \end{align}
    where as defined in assumption \ref{asmp_boundedreward}, there exists $B$ such that $|c_{j}(s,a)-b_j|\leq B$ for all states $s$, actions $a$, and constraints $j \in [m]$.
\end{lemma}
\begin{proof}
    Define once again the vector $g_k := \frac{1}{T_0}\sum_{t=kT_0}^{(k+1)T_0-1} c(s^t,a^t)-b$ as the rollout at epoch $k$, when agents execute $\pi^{k}$. Here $g_k \in \mbR^{m}$, where $m$ is the number of constraints on the system. Then, write the square of the dual iterates' norm at time $k+1$. 
    Using \eqref{eq_lambdaupdate} and the non-expansiveness property of the projection, we get 
    \begin{align} \label{eq_slackness2}
        \|\lambda_{k+1}\|^{2} \leq \|\lambda_k - \eta g_k\|^{2} = \|\lambda_k\|^{2} + \eta^{2} \|g_k\|^{2} - 2\eta g_k^{T}\lambda_k.
    \end{align}
    Using the fact that the difference $|c_{i}(s,a)-b_i| \leq B$, we get 
    \begin{align} \label{eq_slackness3}
        \|\lambda_{k+1}\|^{2} \leq \|\lambda_k\|^{2} + \eta^{2} B^{2} - 2\eta g_k^{T}\lambda_k.
    \end{align}
    Taking conditional expectation with respect to $\sigma$-algebra $\ccalF_k$, and using the fact that $\lambda_k$ is $\ccalF_k$ measurable, we get that 
    \begin{align} \label{eq_slackness4}
        \mbE(\|\lambda_{k+1}\|^{2}|\ccalF_k) \leq \|\lambda_k\|^2 +\eta^2 B^2 -2\eta \mbE(g_k^{T}\lambda_k | \ccalF_k).
    \end{align}
    Now, we shall use \eqref{eq_slackness4} recursively. 
    See that, if we take conditional expectation on both sides with respect to $\ccalF_{k-1}$ for the previous inequality, we get 
    \begin{align} \label{eq_slackness4B}
        \mbE(\|\lambda_{k+1}\|^{2}|\ccalF_{k-1}) \leq \mbE(\|\lambda_k\|^2|\ccalF_{k-1}) +\eta^2 B^2 \nonumber \\ -2\eta \mbE(g_k^{T}\lambda_k | \ccalF_{k-1}).
    \end{align}
    Applying \eqref{eq_slackness4} at epoch $k-1$, we get that $\mbE(\|\lambda_{k}\|^{2}|\ccalF_{k-1}) \leq \|\lambda_{k-1}\|^2 +\eta^2 B^2 -2\eta \mbE(g_{k-1}^{T}\lambda_{k-1} | \ccalF_{k-1})$. Plugging this inequality into \eqref{eq_slackness4B}, we get
    \begin{align} \label{eq_slackness4C}
        \mbE(\|\lambda_{k+1}\|^{2}|\ccalF_{k-1}) \leq \|\lambda_{k-1}\|^2 +2\eta^2 B^2 -2\eta (\mbE(g_k^{T}\lambda_k | \ccalF_{k-1}) \nonumber \\ + \mbE(g_{k-1}^{T}\lambda_{k-1} | \ccalF_{k-1})).
    \end{align}
    Thus, applying \eqref{eq_slackness4} in a recursive fashion, we get that 
    \begin{align} \label{eq_slackness5}
        \mbE(\|\lambda_{k+1}\|^{2}|\ccalF_0) \leq \|\lambda_{0}\|^2 + (k+1)\eta^2 B^2 -2\eta \sum_{l=0}^{k} \mbE(\lambda_l^{T}g_l | \ccalF_0).
    \end{align}
    Thus, using the fact that $\|\lambda_k\| \geq 0$ for all $k$, and rearranging terms in \eqref{eq_slackness5}, we get
    \begin{align}
        \frac{1}{k+1} \sum_{l=0}^{k} \mbE [\lambda_l^T g_l] \leq \frac{\|\lambda_0\|^{2}}{2\eta(k+1)} + \eta \frac{B^2}{2}.
    \end{align}
    Taking limit superior on both sides, we get our final result.
\end{proof}

We now define some terms which will be used to represent the optimal value function for agents against a nonstationary environment comprised of other agents. For every sequence $\{\pi^{k}\}_{k\geq 0}$, we define
\begin{align} \label{def_P_i}
    \ccalP_{i}^{\star} (\pi^{k}_{-i})  :=\sup\limits_{\zeta^{k} \in \Pi^{c}_{i}(\pi^{k}_{-i})} \mbE_{(\zeta^{k},\pi^{k}_{-i})}\left[\sum_{t=kT_0}^{(k+1)T_0-1} r_{i}(s^t,a^t)/T_0\right]
\end{align}
as the optimal value (best-response) for agent $i \in \ccalN$, for epoch $k$, when other agents play $\pi^{k}_{-i}$. 
Further, we define 
\begin{align}
\ccalP_{i}^{\star}:= \liminf_{K\rightarrow \infty} \frac{1}{K} \sum_{k=0}^{K} \ccalP_{i}^{*} (\pi^{k}_{-i}).    
\end{align}
 
We call $\ccalP_{i}^{\star}$ the best-response threshold for agent $i\in \ccalN$ in constrained Markov game $\ccalG$ associated with each sequence $\{\pi^{k}\}_{k\geq 0}$, given other agents play the sequence of policies $\{\pi^{k}_{-i}\}_{k\geq 0}$.  

\begin{lemma}
    Let assumptions \ref{asmp_superfeasible}, \ref{asmp_boundedreward}, \ref{asmp_unbiasedestimate} hold. In addition, assume that the sequence $\{p(\lambda_k | \lambda_0)\}_{k\geq 0}$ is tight, when agents execute the policy sequence  $\{\pi^{k}\}_{k\geq 0}$. Then the sequence of policies forms a constrained nonstationary $\eta\frac{B^2}{2}$-NE for the game $\ccalG$.
\end{lemma}

\begin{proof}
    Let agents follow policy sequence $\{\pi^{k}\}_{k\geq 0}$. Let the policy sequence be realized until epoch $K$. Recall that the generalized dual function for an agent $i\in \ccalN$, Lagrange multiplier $\lambda_k$, and policies of other agents $\pi^{k}_{-i}$, for $k\leq K$ is given by
    \begin{align} \label{eq_optimality1}
        d_i(\lambda_k,\pi^{k}_{-i})= \max_{z_i \in \Pi_i} \left[V_i(z_i,\pi^{k}_{-i}) + \lambda_{k}^{T}(U(z_i,\pi_{-i}^{\lambda_k})-b)\right]
    \end{align}
    Since $\pi^{k}$ is in $NE(\ccalG(\lambda_k))$, we have that $\pi^{k}_i$ is a best-response to $\pi^{k}_{-i}$, which implies that the resultant maximization in \eqref{eq_optimality1} results in
    \begin{align} \label{eq_optimality2}
        d_i(\lambda_k,\pi^{k}_{-i}) = V_i(\pi^{k}) +\lambda_{k}^{T}(U(\pi^{k})-b).
    \end{align}
    
    \noindent \textit{Claim:} We claim that $\ccalP^{\star}_{i}(\pi^{k}_{-i}) \leq d_i(\lambda_k,\pi^{k}_{-i})$. For the sake of contradiction, assume that $\ccalP^{\star}_{i}(\pi^{k}_{-i})> d_i(\lambda_k,\pi^{k}_{-i})$. Using assumption \ref{asmp_unbiasedestimate} (unbiased rollout), we have that  
    \begin{align} \label{eq_optimality2-1}
        \mbE_{(\zeta^{k},\pi^{k}_{-i})}\left[\sum_{t=kT_0}^{(k+1)T_0-1} r_{i}(s^t,a^t)/T_0\right] = V_{i}(\zeta_i, \pi^{k}_{-i})
    \end{align}
    Plugging \eqref{eq_optimality2-1} into \eqref{def_P_i}, we write 
    \begin{align}
        \ccalP^{\star}_{i}(\pi^{k}_{-i})=\sup_{\zeta_i \in \Pi^{c}_{i}(\pi^{k}_{-i})}  V_{i}(\zeta_i, \pi^{k}_{-i}).
    \end{align}
    Since $\ccalP^{\star}_{i}(\pi^{k}_{-i})$ is the maximum value for the constrained problem, the policy $\zeta^{\star}_i=\argmax_{\zeta_i \in \Pi^{c}_{i}(\pi^{k}_{-i})} V_{i}(\zeta_i, \pi^{k}_{-i})$ satisfies 
    $U(\zeta^{\star}_{i}, \pi^{k}_{-i}))-b \geq 0$. Moreover, due to the update dynamics for the Lagrange multiplier in \eqref{eq_lambdaupdate}, $\lambda_k \geq 0$. So, $\lambda_{k}^{T}(U(\zeta^{\star}_{i}, \pi^{k}_{-i}))-b)\geq 0$, which implies that 
    \begin{align}\label{eq_optimality3}
        V_i(\zeta^{\star}_{i},  \pi^{k}_{-i}) + \lambda_{k}^{T}(U(\zeta^{\star}_{i}, \pi^{k}_{-i}))-b) \geq V_i(\zeta^{\star}_{i},  \pi^{k}_{-i})=\ccalP^{\star}_{i}(\pi^{k}_{-i})
    \end{align}
    But, by definition of the generalized dual, $d_i(\lambda_k, \pi^{k}_{-i})\geq V_i(\zeta^{\star}_{i},  \pi^{k}_{-i}) + \lambda_{k}^{T}(U(\zeta^{\star}_{i}, \pi^{k}_{-i}))-b)$, which when reconciled with \eqref{eq_optimality3} creates the contradiction. 

    Thus, from the previous claim, we have $\ccalP^{\star}_{i}(\pi^{k}_{-i}) \leq d_i(\lambda_k,\pi^{k}_{-i})$. Taking average over $K$ epochs, and using \eqref{eq_optimality2}, we have 
    \begin{align} \label{eq_optimality4}
        \frac{1}{K}\sum_{k=0}^{K-1} \left(V_i(\pi^{k}) +\lambda_{k}^{T}(U(\pi^{k})-b)\right) \geq \frac{1}{K}\sum_{k=0}^{K-1} \ccalP^{\star}_{i}(\pi^{k}_{-i}) . 
    \end{align}
    where $\pi_k\in NE(\ccalG(\lambda_k))$.
    Substituting $V_i(\pi^{k})= \mbE \left[ \frac{1}{T_0} \sum_{t=kT_0}^{(k+1)T_0-1} r_{i}(s^t,a^t) \left. \right| \lambda_k  \right]$ and rearranging terms we get 
    \begin{align} \label{eq_optimality5}
        \frac{1}{K} \sum_{k=0}^{K-1} \mbE \left[ \frac{1}{T_0} \sum_{t=kT_0}^{(k+1)T_0-1} r_{i}(s^t,a^t) \left. \right| \lambda_k  \right] \nonumber \\ \geq \frac{1}{K}\sum_{k=0}^{K-1} \ccalP^{\star}_{i}(\pi^{k}_{-i}) - \frac{1}{K} \sum_{k=0}^{K-1} \lambda_{k}^{T}(U(\pi^{k})-b)
    \end{align}
    Now, we take limit on both sides for $K\rightarrow \infty$. We use the result from Lemma \ref{lem_optimality1} to bound $\liminf_{K\rightarrow \infty} \frac{1}{K} \sum_{k=0}^{K-1} \lambda_{k}^{T}(U(\pi^{k})-b) \leq \eta B^2/2$. Recall that $\liminf_{K\rightarrow \infty} \frac{1}{K}\sum_{k=0}^{K-1} \ccalP^{\star}_{i}(\pi^{k}_{-i})= \ccalP^{\star}_i$. Further substituting for $T=KT_0$ on the LHS of \eqref{eq_optimality5}, we get
    \begin{align}
        \liminf_{T\rightarrow \infty} \frac{1}{T} \mbE \sum_{t=0}^{T-1} r_i(s^t, a^t) \geq \ccalP^{\star}_i - \eta B^2 / 2.
    \end{align}
    Since agent $i$ is arbitrary, this shows that the sequence of policies $\{\pi^{k}\}_{k\geq 0}$ comprises a nonstationary $\eta \frac{B^2}{2}$-NE for constrained Markov game $\ccalG$. Moreover, feasibility is guaranteed from proposition \ref{prop_feasibility1}, which in turn implies our desired result.
\end{proof}

\section{Numerical Experiments}
\begin{figure}
    \centering
    \includegraphics[width=0.35\linewidth]{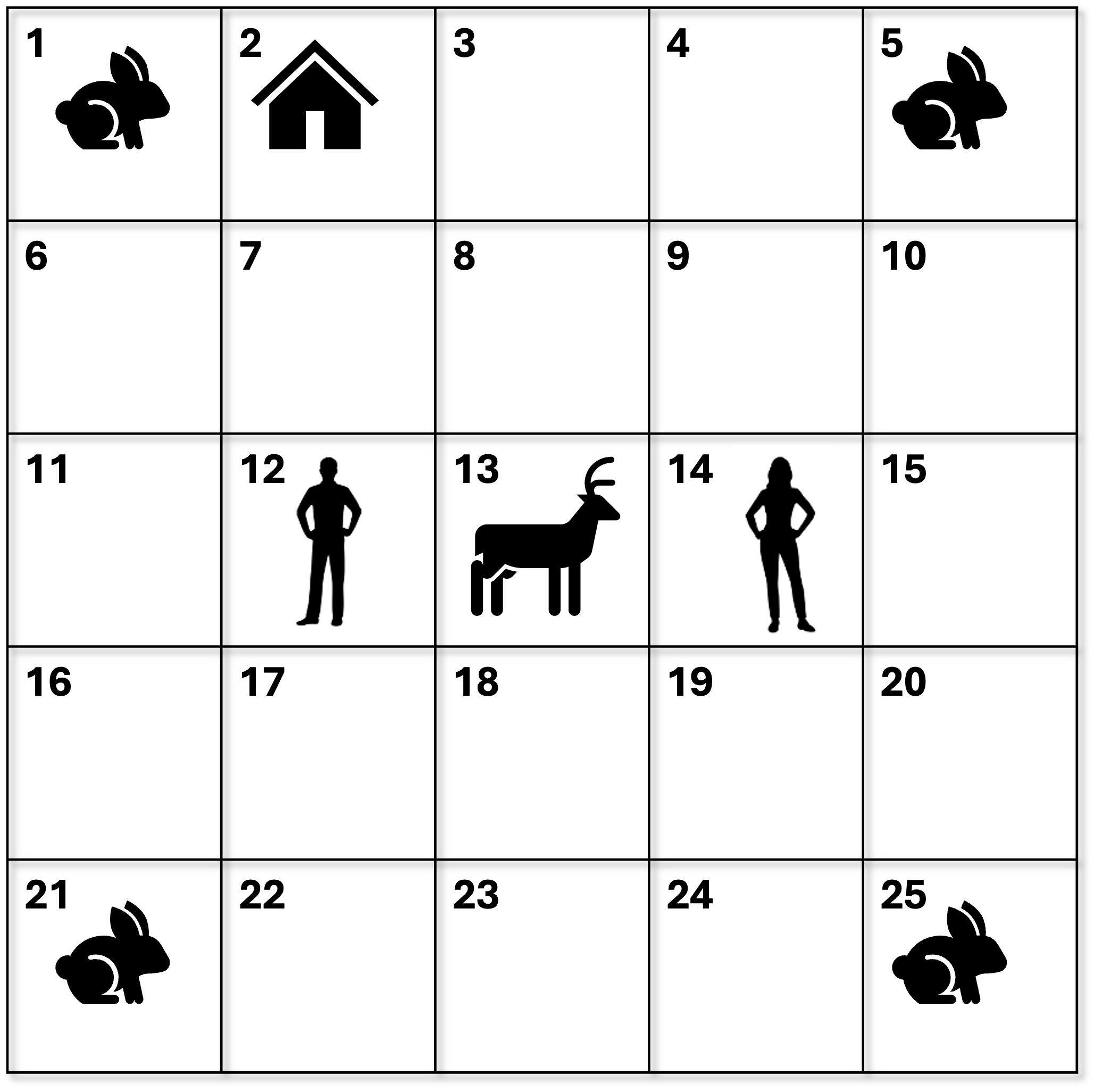}
    \caption{The geometry of the grid for the Stag-Hare-Resting Station game. The position of the players represent the initial state of the game for the trajectory simulations.}
    \label{fig:grid}
\end{figure}

In this example, we study the Stag-Hare-Rest (\textit{SHR}) game. Here two hunters (agents) must coordinate their movements in a grid world to hunt either a stag (which requires cooperation) or hares (individually achievable). Moreover, the agents are required to spend a fraction of their hunting time at a resting station. Each agent can independently choose their next position (local state) by moving to a neighboring grid, but their payoffs depend on their combined position (joint state), and their decision to hunt a stag, or hare, or stay at the resting station. Agents share a common value function (identically interested agents), since the returns are shared among the hunters. Therefore, this setting can be equivalently interpreted as a constrained multiagent MDP.

Additionally, the environment affects agents through the environment's natural transition probabilities, which causes agents to drift away in the grid world according to some state-dependent probability distribution. Agents therefore incur a control cost which quantifies the effort required to alter the environment's natural transitions to align with their navigational strategies. This cost is measured as a KL-Divergence between the stochastic policies of the agents and the natural transition probability distribution, which agents incur at each step along with rewards for hunting hare or stag. 

Let the constrained Markov game example be represented by the tuple $\ccalG^{+}=(\ccalS,\ccalN,\{\ccalA_i, r_i\}_{i\in \ccalN}, P, \{c_{j},b_j\}_{j=1}^{m})$. Then, we have $\ccalN=\{1,2\}$ as the set of agents, and $\ccalS=S_1 \times S_2$ represents the joint positions of the agents in the $5 \times 5$ grid world, where $S_1=S_2=\{1,2,...,25\}$. The location of the hares, the stag and the resting station are as depicted in Figure \ref{fig:grid}. In our implementation, $S_h=\{1,5,21,25\}$ represents the location of hares, $S_s=\{13\}$ is the location of stag, and $S_r=\{2\}$ is the resting station location. Agents have the following actions available to them, namely $\ccalA_i=\{\uparrow, \downarrow,\leftarrow,\rightarrow \}$ for all $i\in \ccalN$. 

The agents receive a reward of $20$ units for hunting stag and $2$ units for hunting hare. Thus, for $s=(s_1,s_2) \in S_1\times S_2$, we have instantaneous reward $r_i(s,a)=20(\mbI_{s_1\in S_r}\mbI_{s_2\in S_r})+2(\mbI_{s_1\in S_h} + \mbI_{s_2 \in S_h})$. Here $\mbI$ stands for the indicator function, and it returns $1$ when the condition (in subscript) is satisfied, and $0$ otherwise. Thus, if both agents are at stag, agents get instantaneous reward of $20$ units, whereas agents incur a reward of $2$ units for each agent at a hare location. Moreover, we have the number of constraints $m=1$ in the game $\ccalG^{+}$. The constraint cost function $c_1$ satisfies $c_1(s,a)=\mbI_{s_1\in S_r}+\mbI_{s_2\in S_r}$. The constraint constant $b_1$ takes values in the set $\{0.25,0.5,0.75\}$ for three different implementations regarding the resting threshold. The transition kernel $P$ responds to agent actions in a deterministic fashion, i.e., if agents select $a_1=\uparrow$ and $a_2=\rightarrow$, the joint state transitions to the new state respecting agent actions, with probability $1$.

Any stationary policy $\pi_i:\ccalS\rightarrow \Delta(\ccalA_i)$ employed by agent $i\in \ccalN$ respects the grid geometry by assigning $0$ probability mass to unallowable actions at state $s\in \ccalS$. For example, an agent situated at the top-right corner of the grid can only move downwards or left. All joint policies in this setup are product policies.
Given $s_i\in S_i$, let $n(s_i)$ represent the adjacent states in the grid world. Let $P^{0}_{i}:S_i\rightarrow \Delta(\ccalA_i)$ represent the natural transition probability distribution for agent $i\in \ccalN$. 
$P^{0}_i$ instructs agent $i$ to retain current state $s_i\in S_i$ with probability $0.9$, and transition to one of the adjacent states $s'\in n(s_i)$ with probability $\frac{1}{|n(s_i)|}$. Then the joint natural transition probability distribution $P^{0}:\ccalS\rightarrow \Delta(\ccalA)$ is given by $P^{0}(s)=P^{0}_{1}(s_1)\times P^{0}_{2}(s_2)$ for $s=(s_1,s_2)\in \ccalS$. Considering the KL-cost, the overall instantaneous reward for agent $i$ for joint state $s\in \ccalS$ and joint action $a\in \ccalA$ is given by $\bar r_i(s,a)=r_i(s,a)+KL(\pi(s)||P^{0}(s))$. In game $\ccalG^{+}$, agents want to maximize $V_{i}^{s}(\pi)=\liminf_{T\rightarrow \infty} \mbE_\pi \frac{1}{T}\sum_{t=0}^{T} \bar r_i(s,a)$ over an infinite horizon of repeated game play, subject to satisfying the system-level constraint $U_1(\pi)=\liminf_{T\rightarrow \infty} \mbE_\pi \frac{1}{T}\sum_{t=0}^{T} c_1(s,a) \geq b_1$. The policy $\pi=\{\pi^{t}\}$ can be nonstationary. For Lagrange multiplier $\lambda$, agents in game $\ccalG^{+}(\lambda)$ have the following Lagrangian value function,  $\ccalL(\pi,\lambda)=V_i(\pi)+\lambda(U_1(\pi)-b_1)$.

As agents simulate Algorithms \ref{alg:dualdes}, \ref{alg:gamedyn}, they move through a sequence of unconstrained Lagrangian games determined by the Lagrange multipliers. 
For any Lagrangian game, we employ an optimistic policy iteration scheme, implemented using temporal-difference based policy evaluation, as the \textit{oracle} to compute the Nash Equilibrium (Assumption \ref{asmp_oracle}). The approach is similar to the one in \cite{nakhleh2024simulation}, \cite{tsitsiklis2002convergence}. Since agents have common value functions, solving for an optimal policy that maximizes the value for an agent in a Lagrangian game also yields a Nash Equilibrium. That is because the unconstrained Lagrangian game coincides with a multiagent MDP in the identical interest cooperative setup. The class of cooperative games can in general be equivalently modeled as optimal control problems by determining the joint value function that agents want to maximize. However, the \textit{oracle} in Assumption \ref{asmp_oracle} is agnostic to how agents reach the stationary NE policy for the associated Lagrangian game. That is, agents are free to employ sophisticated learning mechanisms (such as fictitious play) based on repeated interactions and game-play to reach a stationary NE. Alternatively, a system designer can compute the stationary NE by, for example, solving a constrained MDP in the identical interest setup. Thus, in this example, any other suitable method could be substituted for solving the unconstrained Lagrangian game, such as policy gradients and its variants. 

\begin{figure}
    \centering
    \includegraphics[scale=0.5]{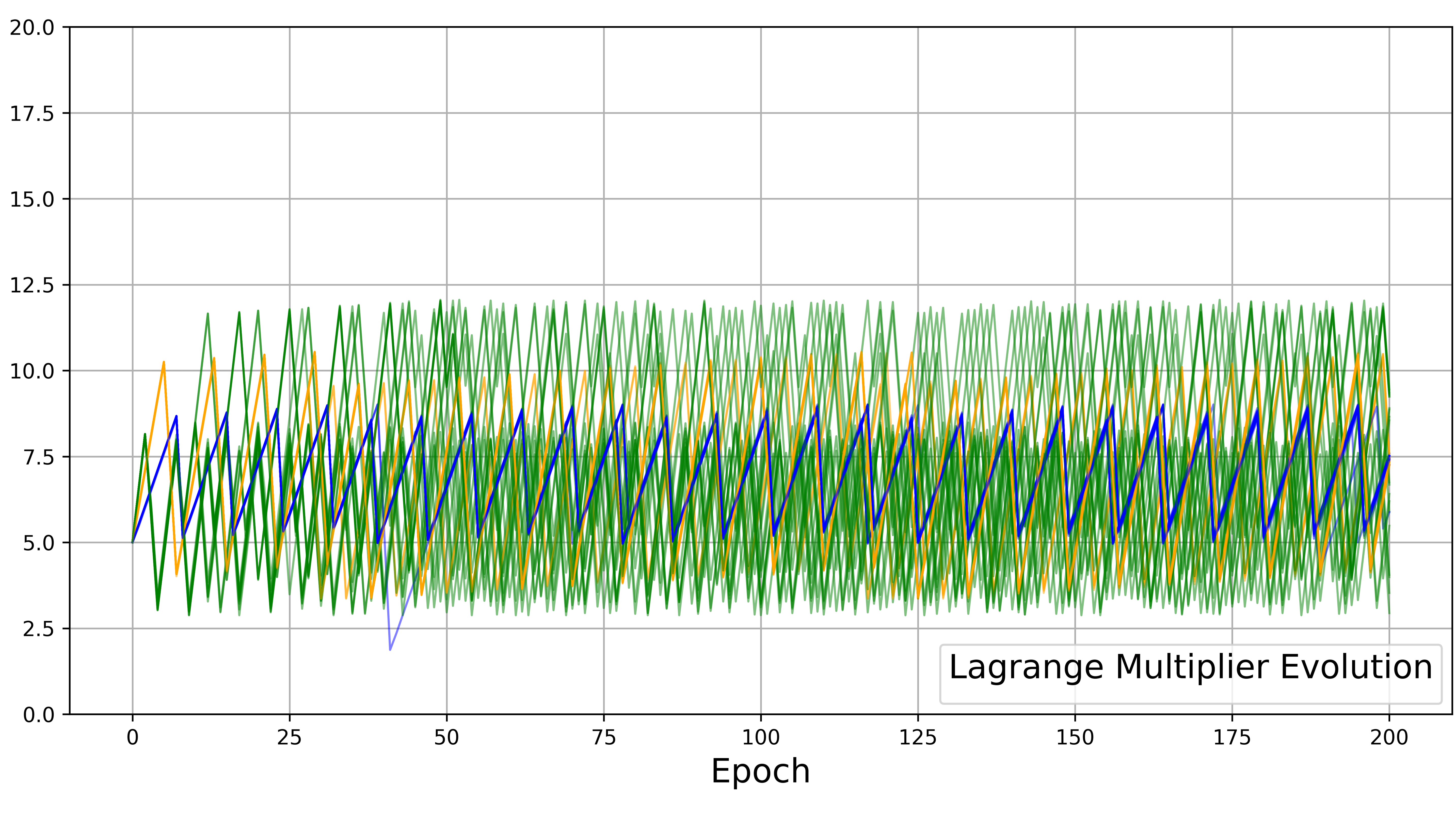}
    \caption{The evolution of the Lagrange multipliers across $K=200$ epochs for multiple episodes sampled from random initial states, with resting threshold levels 0.25 (blue), 0.5 (orange) and 0.75 (green). The multipliers satisfy tightness in Lemma \ref{lem_tightness}.}
    \label{fig:lambda-evol}
\end{figure}
\begin{figure}
    \centering
    \includegraphics[width=0.42\linewidth]{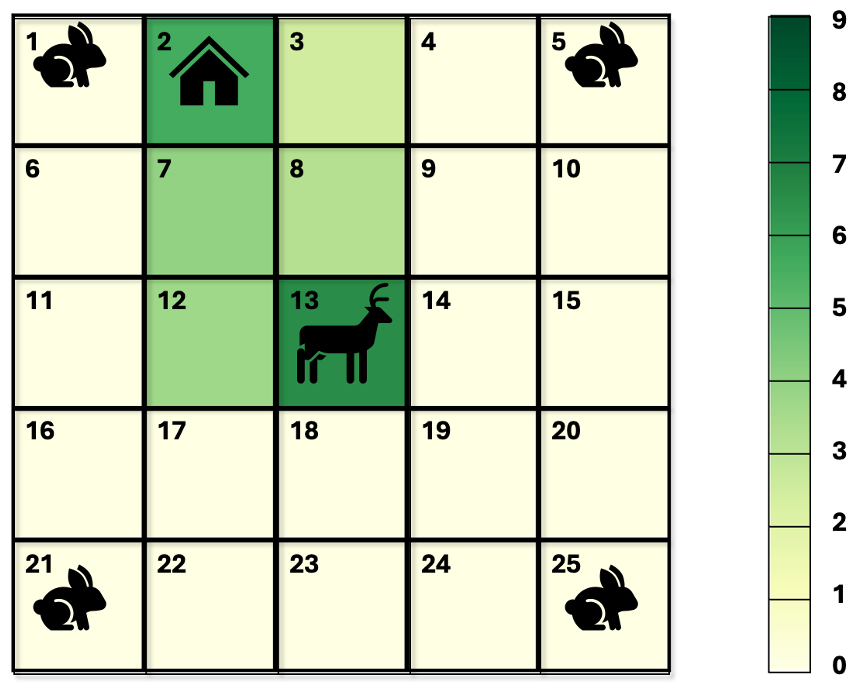}
    \caption{Heat map (on a natural logarithm scale) for the states frequented by agents along an episode. Starting point of agents is $(12,14)$. Agents frequent between hunting stag and resting at the rest-station. They also spend some time traveling between the stag and the resting point.
    Resting threshold level set to $0.50$ for this simulation. Number of epochs simulated is $K=200$ and size of each epoch is $T_0=5$, which implies a total episode length of $t=1000$ time-steps.  }
    \label{fig:Heat Map}
\end{figure}

In Figure \ref{fig:lambda-evol}, we see the oscillation of the Lagrange multipliers simulated according to the dual-descent dynamics in Algorithm \ref{alg:dualdes}. We simulate multiple episodes which start from random initial starting points for the agents, with an initial value of $\lambda_0=5$. 
By simulating Algorithms \ref{alg:dualdes}, \ref{alg:gamedyn}, we construct an episode of agents playing $\{\pi^{k}\}_{k\geq0}$, where for each epoch $k\leq K$, agents play policy $\pi^{k}\in NE(\ccalG(\lambda_k))$. 
Starting at $\lambda_0=5$, agents initially accrue constraint violations for not satisfying the resting requirement. At smaller values of the multiplier, agents play policies that prioritize hunting stag over resting. However, this causes the Lagrange multipliers to go up, following the dual-descent dynamics. At larger values of the Lagrange multiplier, the agents start frequenting the resting station, by playing stationary policies which prioritize resting. This leads to the oscillation we see in Figure \ref{fig:lambda-evol}. We see another demonstration of this oscillatory phenomena in Figure \ref{fig:Heat Map}. The heat-map plots the occupancy measure resulting from running an episode according to Algorithm \ref{alg:gamedyn}. The map strongly suggests that agents spend most of their time hunting stag, followed by resting, followed by traveling between the two locations. 




\begin{figure} 
    \centering
    \includegraphics[scale=0.5]{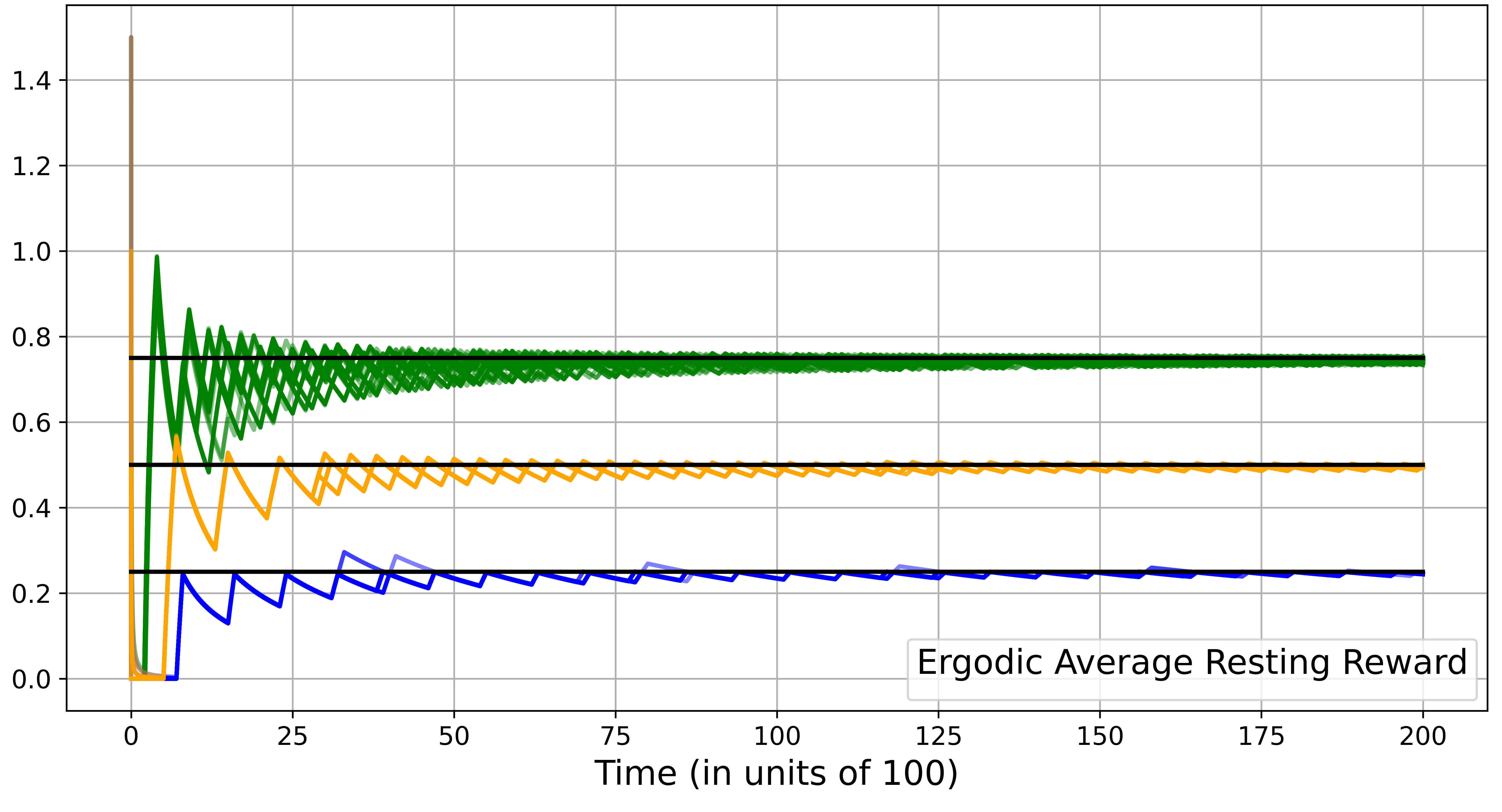}
    \caption{Constraint satisfaction by the simulated episodes as the resting station rewards (total time spent at resting station by both the agents) collapses on the resting time threshold (threshold values 0.25, 0.50 and 0.75 in blue, orange and green respectively). Episodes simulated from random initial states. Number of epochs $K=200$, size of epoch $T_0=100$.}
    \label{fig:const_sat}
\end{figure}

\begin{figure}
    \centering
    \includegraphics[scale=0.5]{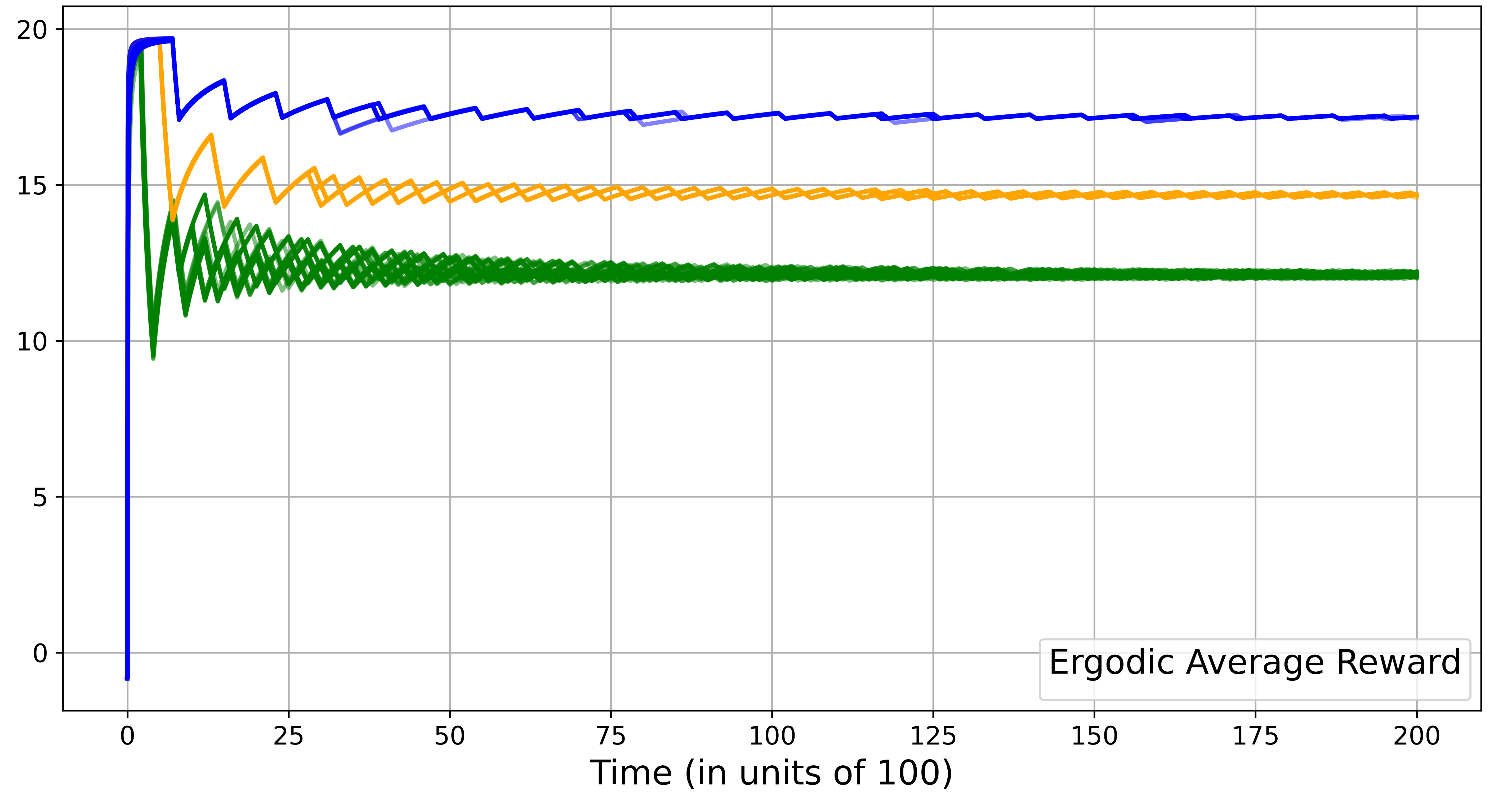}
    \caption{Time average cumulative value episodes converge for agents, constituting an $\epsilon$-NE. The blue, orange and green curves represent resting thresholds of 0.25, 0.50 and 0.75 respectively. Episodes simulated from random initial states. Number of epochs $K=200$, size of epoch $T_0=100$.}
    \label{fig:val_con}
\end{figure}

These oscillations between policies that prioritize hunting stag versus policies that prioritize resting provide the foundation for the constraint satisfaction results we observe in Figure \ref{fig:const_sat}. Starting from random initial states, the green, orange and blue curves correspond to time average resting rewards as a function of time (on the horizontal axis), obtained as a result of the execution of Algorithm \ref{alg:gamedyn}. The resting reward curves oscillate around their individual resting thresholds ($0.75$ for green, $0.5$ for orange, $0.25$ for blue), and eventually collapse on these thresholds, proving constraint satisfaction with no slack. Moreover, the corresponding value curves (see Figure \ref{fig:val_con}) converge too, with the values converging to a higher magnitude for the case when the resting threshold requirement is lower. This indicates, that when the resting requirement is lower, agents spend more time hunting stag, and thus accrue a higher ergodic average reward. 

The tight constraints ensure that agents maximize hunting rewards across an episode while satisfying the constraints. While our theoretical guarantees do not feature optimal value maximization, as optimality is ill-defined for a constrained Markov game in general, we still achieve optimality for the constrained multiagent MDP example under consideration. That is, in this case where agents are identically interested, as in the \textit{SHR} example, our results in Figures \ref{fig:const_sat},\ref{fig:val_con} suggest that the sequence of policies $\{\pi^k\}_{k\geq 0}$ constitutes an approximately optimal nonstationary policy. This results from the fact that for each epoch $k$, our optimistic policy iteration algorithm yields a stationary policy $\pi^{k}$ that approximately maximizes the common value function for the agents. This is also evident from the heat-map for the states visited along an episode (Figure \ref{fig:Heat Map}), as agents ignore individually hunting hares for the sake of switching between hunting stag and resting at the rest-stop.

\section{Conclusion}
In this paper, we present the formulation of a Lagrangian game, given an associated constrained Markov game. We show that solving a sequence of carefully generated Lagrangian games produced by continuous execution of the rollout dual descent (Algorithm 1) and game dynamics (Algorithm 2) produces cost trajectories which are feasible almost surely. Moreover, the associated sequence of policies, which are individually stationary solutions for the series of Lagrangian games, yields a nonstationary Nash solution for the original constrained Markov game. Thus, our results produce a systematic way to exploit primal dual methods in general constrained Markov games, marking a first step in the literature. We provide numerical simulation results for a identical interest constrained Markov game to back up our theoretical results. Our framework is novel in that it reduces the problem of solving a constrained Markov game to a sequence of unconstrained Markov games, significantly broadening the practical applicability of constrained game solutions without relying on sophisticated optimization techniques for handling constraints directly.
\bibliographystyle{IEEEtran}
\bibliography{Main}
\end{document}